\documentclass[11pt]{article}
\usepackage{amssymb}
\usepackage[dvips]{epsfig}
\usepackage{graphicx}

\newtheorem{Lemma}{Lemma}[section]
\newtheorem{Theorem}{Theorem}
\newtheorem{Proposition}[Lemma]{Proposition}

\newtheorem{Remark}[Lemma]{Remark}

\makeatletter\@addtoreset{figure}{section}\makeatother

\makeatletter\@addtoreset{equation}{section}\makeatother

\newcommand{\C}{\mathbb{C}}

\newcommand{\R}{\mathbb{R}}
\newcommand{\Z}{\mathbb{Z}}
\newcommand{\T}{\mathbb{T}}

\def\Re{\mathop{\mathrm{Re}}}
\def\Im{\mathop{\mathrm{Im}}}

\newcommand{\rmO}{\mathcal{O}}

\newcommand{\rmd}{\mathrm{d}}
\newcommand{\rme}{\mathrm{e}}
\newcommand{\rmi}{\mathrm{i}}

\newcommand{\per}{\mathrm{per}}
\newcommand{\loc}{\mathrm{loc}}
\newcommand{\sign}{\mathrm{sign}}
\newcommand{\linspan}{\mathrm{span}}
\newcommand{\Det}{\mathrm{det}}

\newcommand{\dd}{\,\mathrm{d}}
\renewcommand{\AA}{\mathcal{A}}
\newcommand{\BB}{\mathcal{B}}
\newcommand{\EE}{\mathcal{E}}
\newcommand{\MM}{\mathcal{M}}
\newcommand{\HH}{\mathcal{H}}
\newcommand{\KK}{\mathcal{K}}
\newcommand{\TT}{\mathcal{T}}
\newcommand{\RR}{\mathcal{R}}
\newcommand{\NN}{\mathcal{N}}
\newcommand{\bA}{\mathbf{A}}
\newcommand{\bB}{\mathbf{B}}
\newcommand{\bC}{\mathbf{C}}
\newcommand{\1}{\mathbf{1}}

\renewcommand{\phi}{\varphi}

\newcommand{\reff}[1]{(\ref{#1})}

\def\eqdef{\buildrel\hbox{\small{def}}\over =}
\def\build#1_#2^#3{\mathrel{
  \mathop{\kern 0pt#1}\limits_{#2}^{#3}}}
\def\QED{\mbox{}\hfill$\Box$}

\newdimen\texpscorrection
\newdimen\figcenter
\def\figurewithtex #1 #2 #3 #4 #5\cr{\null
  {\goodbreak\figcenter=\hsize\relax
  \advance\figcenter by -#4truecm
  \divide\figcenter by 2
  \begin{figure}[hbt]
  \vskip #3truecm\noindent\hskip\figcenter
  \includegraphics{#1}{\hskip\texpscorrection\input #2 }
  \vskip 0.8truecm{\baselineskip=0.8\baselineskip
  \noindent \vbox{\noindent {\footnotesize #5}}\par}
  \end{figure}}}
\def\point#1 #2 #3 {\rlap{\kern #1 truecm
\raise #2 truecm \hbox{#3}}}

\begin{document}
\begin{center}

\vspace*{4mm}

{\LARGE Orbital stability of periodic waves for the nonlinear
\\[1ex] Schr\"odinger equation}
\\[8mm]
{\large\bf Thierry Gallay}
\\
Institut Fourier\\
Universit{\'e} de Grenoble I\\
B.P. 74 \\
38402 Saint-Martin-d'H{\`e}res, France
\\[4mm]
{\large\bf Mariana  H\u{a}r\u{a}gu\c{s}}
\\
D\'epartement de Math\'ematiques\\
Universit\'e de Franche-Comt\'e\\
16 route de Gray\\
25030 Besan\c{c}on, France
\end{center}

\vspace*{1cm}

\begin{abstract}
The nonlinear Schr\"odinger equation has several families of
quasi-periodic travelling waves, each of which can be parametrized
up to symmetries by two real numbers: the period of the modulus of
the wave profile, and the variation of its phase over a period
(Floquet exponent). In the defocusing case, we show that these
travelling waves are orbitally stable within the class of solutions
having the same period and the same Floquet exponent. This
generalizes a previous work \cite{GH1} where only small amplitude
solutions were considered. A similar result is obtained in the
focusing case, under a non-degeneracy condition which can be checked
numerically. The proof relies on the general approach to orbital
stability as developed by Grillakis, Shatah, and Strauss
\cite{GSS1,GSS2}, and requires a detailed analysis of the
Hamiltonian system satisfied by the wave profile.
\end{abstract}

\thispagestyle{empty}

\vfill {\bf Running head:} Periodic waves in the NLS equation

{\bf Corresponding author:} Thierry Gallay,
{\tt Thierry.Gallay@ujf-grenoble.fr}

{\bf Keywords:} Nonlinear Schr\"odinger equation, periodic waves,
orbital stability

\newpage

\section{Introduction}
\label{s:1}

This paper is devoted to the stability analysis of the quasi-periodic
travelling wave solutions of the cubic nonlinear Schr\"odinger
(NLS) equation
\begin{equation}\label{e:nls}
  \rmi U_t(x,t) + U_{xx}(x,t) + \gamma |U(x,t)|^2 U(x,t) \,=\, 0~,
  \quad x \in \R~, \quad t \in \R~,
\end{equation}
where $\gamma \in \{-1;1\}$ and $U(x,t) \in \C$.  Eq.\reff{e:nls}
is a universal envelope equation describing the propagation of
weakly nonlinear waves in dispersive media (see \cite{Su} for a
comprehensive introduction).  The nonlinearity in \reff{e:nls} is
``attractive'' if $\gamma = +1$ (focusing case) and ``repulsive''
if $\gamma = -1$ (defocusing case).  In both cases Eq.\reff{e:nls}
has a family of quasi-periodic travelling waves of the form
\begin{equation}\label{e:qpdef}
  U(x,t) \,=\, \rme^{\rmi (p x - \omega t)}\,V(x-ct)~,
  \quad x \in \R~, \quad t \in \R~,
\end{equation}
where $p,\omega,c$ are real parameters and $V : \R \to \C$
is a periodic function. The simplest elements of this family are
the {\it plane waves}, for which $V$ is identically constant and
$p,\omega$ satisfy the dispersion relation $p^2 = \omega + \gamma
|V|^2$. It is well-known and easy to verify that the plane waves
are dynamically stable in the defocusing case, and unstable (if
$V \neq 0$) in the focusing case \cite{Z,GH1}. We shall therefore
concentrate on the less explored situation where $V$ is
a nontrivial periodic function. In that case, we shall refer
to \reff{e:qpdef} as a {\it periodic wave}, although $U(x,t)$ is
in general a quasi-periodic function of both $x$ and $t$.

The number of parameters in \reff{e:qpdef} can be reduced if we use
the symmetries of Eq.\reff{e:nls}. We recall that the NLS equation is
invariant under the following transformations:
\begin{enumerate}
  \item $U(x,t) \mapsto U(x,t)\,\rme^{\rmi\phi}$, $\phi\in\R$
  (Phase invariance);
  \item $U(x,t) \mapsto U(x+\xi,t)$, $\xi\in\R$ (Translation
    invariance);
  \item $U(x,t)\mapsto \rme^{-\rmi\big(\frac v2x + \frac{v^2}4t\big)}
    U(x+vt,t)$, $v\in\R$ (Galilean invariance);
  \item $U(x,t) \mapsto \lambda U(\lambda x,\lambda^2t)$, $\lambda > 0$
    (Dilation invariance).
\end{enumerate}
If $U(x,t)$ is a periodic wave as in \reff{e:qpdef}, we can use
the Galilean invariance to transform it into a solution of the
same form with $c = 0$. Then, using the dilation invariance, we
can further assume that $\omega \in \{-1;0;1\}$. It follows that
$U(x,t) = \rme^{-\rmi \omega t} W(x)$, where $W(x) = \rme^{\rmi p x}
V(x)$ is a solution of the ordinary differential equation
\begin{equation}\label{e:snls}
  W_{xx}(x) + \omega W(x) + \gamma |W(x)|^2 W(x) \,=\, 0~, \quad
  x \in \R~.
\end{equation}

The bounded solutions of \reff{e:snls} are completely classified
for all values of the parameters $\omega, \gamma$. The simplest
ones are the plane waves $W(x) = A\,\rme^{\rmi px}$ where $p \in
\R$, $A \in \C$ and $p^2 = \omega + \gamma |A|^2$. The periodic
waves correspond to quasi-periodic solutions of \reff{e:snls} of
the form $W(x) = r(x)\,\rme^{\rmi \phi(x)}$, where $r,\phi$ are
real functions with the property that $r$ and $\phi_x$ are
periodic with the same period. It turns out that Eq.\reff{e:snls}
has a four-parameter family of such solutions, both in the
focusing and in the defocusing case (see Sections~\ref{s:ex} and
\ref{s:fe} below). Actually, if $\gamma = -1$, we must assume that
$\omega = 1$ otherwise \reff{e:snls} has no nontrivial bounded
solutions; if $\gamma = +1$, Eq.\reff{e:snls} has quasi-periodic
solutions for all values of $\omega$, but we shall only consider
the generic cases $\omega = \pm 1$. If $\gamma \omega < 0$, in
addition to plane waves and periodic waves, there exist pulse-like
solutions of \reff{e:snls} which are homoclinic as $x \to
\pm\infty$ to a plane wave or to the zero solution. The most
famous one (if $\gamma = 1$ and $\omega = -1$) is the ground state
$W(x) = \sqrt{2}/\cosh(x)$ which corresponds to the solitary wave
of the focusing NLS equation.

In contrast to the plane waves or the solitary waves which have been
extensively studied \cite{CL,W,Z}, relatively little seems to be known
about the stability of periodic waves. {\em Spectral stability} with
respect to long-wave disturbances has been examined by Rowlands
\cite{R}, who showed that periodic waves with real-valued profile are
unstable in the focusing case and stable (at least in the long-wave
regime) in the defocusing case. In the latter case, spectral stability
of the small amplitude periodic waves has been rigorously established
in \cite{GH1} using Bloch-wave analysis. Similar results were also
obtained for certain NLS-type equations with spatially periodic
potentials \cite{BR2,OSY}. As for the {\em nonlinear stability}, the
only result we are aware of is due to Angulo \cite{Ang}, who proved
very recently that the family of dnoidal waves of the focusing NLS
equation is orbitally stable with respect to perturbations which have
the same period as the wave itself. We recall that the periodic waves
of NLS with real-valued profile are called ``cnoidal waves'' when they 
have zero average over a period (like the Jacobian elliptic function 
$cn$), and ``dnoidal waves'' when they have nonzero average (like the 
elliptic function $dn$). 

In this paper, we study the nonlinear stability of  all periodic
waves of \reff{e:nls}, but we restrict ourselves to a specific
class of perturbations which we now describe. Any quasi-periodic
solution of \reff{e:snls} can be written in the form
\begin{equation}\label{e:WQ}
  W(x) \,=\, \rme^{\rmi p x}\,Q_\per(2kx)~, \quad x \in \R~,
\end{equation}
where $p \in \R$, $k > 0$, and $Q_\per : \R \to \C$ is $2\pi$-periodic.
Here $k = \pi/T$, where $T > 0$ is the minimal period of $|W|$.
The representation \reff{e:WQ} is not unique, since we can add to
$p$ any integer multiple of $2k$ (and modify the periodic function
$Q_\per$ accordingly), but the Floquet multiplier $\rme^{\rmi p T}$ is
uniquely defined. Our purpose is to show that the periodic wave
$U(x,t) =\rme^{\rmi (px-\omega t)} Q_\per(2kx)$ of \reff{e:nls} is
stable within the class of solutions which have the {\em same
period} $T = \pi/k$ and the {\em same Floquet multiplier}
$\rme^{\rmi p T}$. In other words, we restrict ourselves to
solutions of \reff{e:nls} of the form $\rme^{\rmi (px-\omega t)}
Q(2kx,t)$, where $Q(\cdot,t)$ lies in the function space
\[
  X \,=\, H^1_\per([0,2\pi],\C) \,=\, \Bigl\{Q \in H^1_\loc(\R,\C)
  \,\Big|\, Q(z) = Q(z+2\pi) \hbox{ for all } z \in \R\Bigr\}~.
\]
The advantage of this restricted setting is that nonlinear
stability can be established by the standard variational method
which has been developed originally to prove the orbital stability
of solitary waves \cite{Be,Bo,W} (see also \cite{Ang,ABS}).
However, the obvious drawback of this approach is that it does not
give any information on the stability of the periodic waves with
respect to non-periodic perturbations, a difficult question which
remains essentially open.

With this perspective in mind, we shall put the emphasis on the
{\em defocusing case} $\gamma = -1$, because we know from \cite{R}
that the periodic waves will be unstable in the focusing case if
non-periodic perturbations are to be allowed.  Our main result
can be stated as follows:

\begin{Theorem}\label{th:orbit}
{\rm (Orbital stability of periodic waves in the defocusing case)}\\
Let $\gamma=-1$, $\omega = 1$, and assume that $W(x) = \rme^{\rmi
px}\,Q_\per(2kx)$ is a solution of \reff{e:snls} with $p \in \R$, $k
> 0$, and $Q_\per \in X$, as in \reff{e:WQ}.  Then there exist $C_0 >
0$ and $\epsilon_0 > 0$ such that, for all $R \in X$ with $\|R\|_X \le
\epsilon_0$, the solution $U(x,t) = \rme^{\rmi (px-\omega t)}Q(2kx,t)$
of the NLS equation \reff{e:nls} with initial data $U(x,0) =
\rme^{\rmi px} (Q_\per(2kx) + R(2kx))$ satisfies, for all $t \in \R$,
\begin{equation}\label{e:orbit}
  \inf_{\phi,\xi \in [0,2\pi]} \|Q(\cdot,t) - \rme^{\rmi \phi}
  Q_\per(\cdot - \xi)\|_X \,\le\, C_0 \|R\|_X~.
\end{equation}
\end{Theorem}

This result is known to hold for small amplitude periodic waves
\cite{GH1}, in which case the constants $C_0$, $\epsilon_0$ do not
depend on the wave profile $Q_\per$. Here we remove the smallness
assumption, but our argument relies in part on the calculations
made in \cite{GH1} (see Lemma~\ref{th:smallwaves} below).

\medskip

\noindent
\textbf{Remarks}\\
{\bf 1.} Theorem~\ref{th:orbit} includes the situation where
$|W|$ is constant, in which case $W$ is a {\em plane wave}
rather than a periodic wave. Since stability is well-known for
plane waves \cite{Z}, we shall assume henceforth that $|W|$
is a nontrivial periodic function. In such case, we emphasize 
that the wavenumber $k$ which appears in \reff{e:WQ} is always
given by $k = \pi/T$, where $T > 0$ is the {\em minimal} period
of $|W|$. This is very important because our approach does not 
allow to prove the stability of periodic waves with respect 
to perturbations whose period is an integer multiple of $T$.  
\\[1mm]
{\bf 2.} It is interesting to see what Theorem~\ref{th:orbit} means in
the particular case of cnoidal waves. For such waves we have $W(x) =
\rme^{\rmi px} Q_\per(2kx)$ where $p=k=\pi/T$ and $T>\pi$ is the
minimal period of $|W|$. The Floquet multiplier $\rme^{\rmi pT}$ is
therefore equal to $-1$, so that $W(x+T)=-W(x)$, for all $x\in\R$. In
particular, $W$ is periodic with (minimal) period $L=2T$.
Theorem~\ref{th:orbit} then shows that the $L$-periodic cnoidal wave
$U(x,t) = \rme^{-\rmi t}W(x)$ is orbitally stable with respect to
$L$-periodic perturbations $\widetilde W$ {\em provided} that
$\widetilde W(x+L/2) = - \widetilde W(x)$ for all $x\in\R$. As
explained in \cite{Ang}, without this additional assumption the
classical approach does not allow to prove the orbital stability of
cnoidal waves with respect to perturbations which have the same period
as the wave itself.

\medskip

The rest of the paper is organized as follows. In Section~\ref{s:ex}
we recall the classification of the bounded solutions of \reff{e:snls}
with $\gamma = -1$ and $\omega = 1$. These solutions can be
interpreted as the trajectories of an integrable Hamiltonian system
with two degrees of freedhom, which is proved to be non-degenerate in
the sense of KAM theory. We thereby answer a question raised by
Bridges and Rowlands \cite{BR1} in connexion with the stability of the
quasi-periodic solutions of the Ginzburg-Landau equation.

Section~\ref{s:os} is devoted to the proof of Theorem~\ref{th:orbit}.
As indicated above, we follow the general approach of Grillakis,
Shatah, and Strauss. The main difficulty is to verify the assumptions
of the stability theorem in \cite{GSS2}. We first check that the
second variation of the energy functional at the periodic wave has
exactly one negative eigenvalue. This result has been established for
small waves in \cite{GH1}, and a continuity argument allows to extend
it to periodic waves of arbitrary size. We next consider the {\em
structure function} (which is called ``$d(\omega)$'' in \cite{GSS2})
and show, by a direct calculation, that its Hessian matrix has a
negative determinant. Both properties together imply orbital
stability.

Finally, in Section~\ref{s:fe}, we extend our results to the
focusing NLS equation. The situation is more complicated here,
because we have families of periodic waves for all values of
$\omega \in \{-1;0;1\}$. By comparison with the spatially
homogeneous rotating wave $U(x,t) = \rme^{\rmi t}$, these periodic
waves may be called ``counter-rotating'' if $\omega = +1$,
``standing'' if $\omega= 0$, and ``corotating'' if $\omega = -1$.
As was already mentioned, we restrict ourselves to the generic
cases $\omega = \pm 1$.  In both situations, we show that the
Hamiltonian system corresponding to \reff{e:snls} is
non-degenerate in the sense of KAM, and we deduce as in
Section~\ref{s:os} that the second variation of the energy
functional has exactly one negative eigenvalue. It follows that
orbital stability holds provided the Hessian matrix of the
structure function has a negative determinant (see
Propositions~\ref{th:orbit1} and \ref{th:orbit2} for precise
statements). Unlike in the defocusing case, we do not give a
complete proof of this property, because the calculations are
excessively complicated. But the determinant is easy to evaluate
numerically (or even analytically in various parameter regimes),
and it appears to be negative for all periodic waves. Assuming
this to be true, we conclude that the analogue
Theorem~\ref{th:orbit} holds in the focusing case too. Thus, in
contrast to what happens when long-wave disturbances are
considered \cite{R}, there is apparently no difference between the
focusing and the defocusing case as far as periodic perturbations
are considered.

\medskip\noindent{\bf Acknowledgements.} The authors are indebted
to A. De Bouard and L. Di Menza for fruitful discussions. This
work was partially supported by the French Ministry of Research
through grant ACI JC 1039.

\section{Properties of the periodic waves}
\label{s:ex}

In this section, we study the bounded solutions of the stationary
Ginzburg-Landau equation
\begin{equation}\label{e:sgl}
  W_{xx}(x) + W(x) -|W(x)|^2 W(x) \,=\, 0~,
\end{equation}
where $W : \R \to \C$. If we interprete the spatial variable
$x \in \R$ as a ``time'', Eq.\reff{e:sgl} becomes an integrable
Hamiltonian dynamical system with two degrees of freedhom. The
conserved quantities are the ``angular momentum'' $J$ and
the ``energy'' $E$:
\begin{equation}\label{e:JEdef}
  J \,=\, \Im(\overline{W}W_x)~, \qquad
  E \,=\, \frac12|W_x|^2 + \frac12|W|^2 - \frac14|W|^4~.
\end{equation}
If $W$ is a solution of \reff{e:sgl} with $J \neq 0$, then $W(x)
\neq 0$ for all $x \in \R$, so that we can introduce the polar
coordinates $W(x) = r(x)\,\rme^{\rmi \phi(x)}$. The invariants
then become
\[
  J \,=\, r^2 \phi_x~, \qquad E \,=\,\frac{r_x^2}2 +
  \frac{J^2}{2r^2} + \frac{r^2}2 - \frac{r^4}4~.
\]

Let $D \subset \R^2$ be the open region defined by
\begin{equation}\label{e:Ddef}
  D \,=\, \Bigl\{(J,E)\in\R^2 \,\Big|\, J^2 < 4/27\,,~E_-(J) < E
  < E_+(J)\Bigr\}~,
\end{equation}
where the functions $E_-,E_+$ are defined in \reff{e:Epmdef} below
(see also Fig.~1). Our starting point is the following well-known
classification \cite{BR1,DGJ,Ga}:

\figurewithtex fig1.ps Fig1.tex 5.300 11.000
{\bf Fig.~1:} The region $D$ in the parameter space for which
Eq.\reff{e:sgl} has bounded solutions.\cr

\begin{Proposition}\label{th:glstat}
If $W : \R \to \C$ is a bounded solution of \reff{e:sgl}, the
corresponding invariants $(J,E)$ belong to the closure
$\overline{D}$ of $D$ and the following properties hold:
\begin{enumerate}
\item If $(J,E) \in D$ and $J \neq 0$, then the modulus $|W|$ and
the derivative of the phase of $W$ are periodic with the same
minimal period $T(J,E) > 0$. If $(J,E) \in D$ and  $J=0$, then $W$
is real-valued, up to a phase factor, and periodic with period
$2T(0,E)$.

\item If $(J,E) \in
\partial D$, then either $W(x) = W_{q,\phi}(x) \equiv
(1-q^2)^{1/2} \rme^{\rmi (qx +\phi)}$ for some $q \in [-1,1]$ and
some $\phi \in [0,2\pi]$, or $W$ is a homoclinic orbit connecting
$W_{q,\phi_-}$ at $x = -\infty$ to $W_{q,\phi_+}$ at $x =
+\infty$, for some $q^2 < 1/3$ and some $\phi_-,\phi_+ \in
[0,2\pi]$.
\end{enumerate}
\end{Proposition}

\noindent{\bf Proof.}
Although the arguments are standard, we give a complete proof of
Proposition~\ref{th:glstat} because it will serve as a basis
for all subsequent developments in this section.
Let $W : \R \to \C$ be a bounded solution of \reff{e:sgl}, and
assume first that $J \neq 0$. Then $W(x) = r(x)\,\rme^{\rmi \phi(x)}$
and $E = \frac12 r_x^2 + V_J(r)$, where $V_J$ is the
``effective potential''
\[
  V_J(r) \,=\, \frac{J^2}{2r^2} + \frac{r^2}2 - \frac{r^4}4~, \quad
  r > 0~.
\]
If $J^2 > 4/27$, then $V_J'(r) < 0$ for all $r > 0$, hence
\reff{e:sgl} has no bounded solution in that case. Thus we must
have $J^2 \le 4/27$. If $0 < J^2 < 4/27$, we can parametrize $J$
in a unique way as
\[
  J \,=\, q(1-q^2) \,=\, Q(1-Q^2)~, \quad \hbox{where} \quad
  0 < q^2 < 1/3 < Q^2 < 1~.
\]
Then $y^3 - y + J = (y-q)(y-Q)(y+q+Q)$, hence in particular $J =
qQ(q+Q)$. With this parametrization, it is easy to check that
$V_J(r)$ has a (unique) local minimum at $r_Q = (1-Q^2)^{1/2}$ and
a (unique) local maximum at $r_q = (1-q^2)^{1/2}$ (see Fig.~2). We
define
\begin{eqnarray}\label{e:Epmdef}
  E_-(J) &=& V_J(\sqrt{1-Q^2}) \,=\, \frac14(1-Q^2)(1+3Q^2)~,\\ \nonumber
  E_+(J) &=& V_J(\sqrt{1-q^2}) \,=\, \frac14(1-q^2)(1+3q^2)~.
\end{eqnarray}
Since $W$ is a bounded solution of \reff{e:sgl}, we necessarily have
$E_-(J) \le E \le E_+(J)$. This gives three possibilities:

\begin{enumerate}
\item If $E = E_-(J)$, then (up to a global phase factor)
$W(x) = (1-Q^2)^{1/2}\,\rme^{\rmi Qx}$, i.e. $W$ is a periodic solution
with constant modulus.
\item If $E = E_+(J)$, then $|W|$ is either constant or homoclinic
to $r_q$ as $x \to \pm\infty$. In the first case, $W(x) =
(1-q^2)^{1/2}\,\rme^{\rmi qx}$ (up to a phase factor). In the second
case, up to a translation and a phase factor, we have the explicit
formula
\[
  W(x) \,=\, \Bigl(2(q^2 + c^2 \tanh^2(cx))\Bigr)^{1/2}
  \,\rme^{\rmi qx + \rmi \arctan (\frac{c}{q}\tanh(cx))}~, \quad
  \hbox{where } c \,=\, \sqrt{\frac{1-3q^2}{2}}~.
\]
\item If $E_-(J) < E < E_+(J)$, the modulus $r = |W|$ and
the phase derivative $\phi_x = \Im(W_x/W)$ are periodic with the same
period. If we denote by $r_1 < r_2 < r_3$ the three positive
roots of $E - V_J(r)$ as in Fig.~2, this (minimal) period is
\begin{equation}\label{e:Tdef}
  T(J,E) \,=\, 2\int_{r_1(J,E)}^{r_2(J,E)} \frac{\dd r}
  {\sqrt{2(E - V_J(r))}}~.
\end{equation}
Another important quantity is the increment of the phase
$\phi$ over a period of the modulus, namely
\begin{equation}\label{e:Phidef}
  \Phi(J,E) \,=\, 2\int_{r_1(J,E)}^{r_2(J,E)} \frac{J}{r^2}\frac{\dd r}
  {\sqrt{2(E - V_J(r))}}~.
\end{equation}
Since in general $\Phi(J,E)$ is not a rational multiple of $\pi$,
the solution $W(x) =  r(x)\,\rme^{\rmi \phi(x)}$ of \reff{e:sgl}
is not periodic, but only quasi-periodic.
\end{enumerate}

\figurewithtex 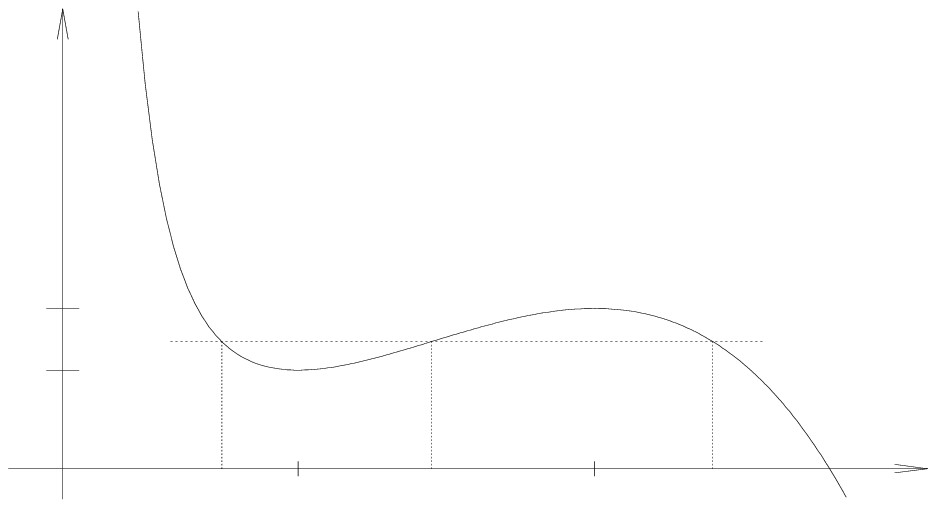 Fig2.tex 5.500 10.000
{\bf Fig.~2:} When $0 < J^2 < 4/27$, the effective potential $V_J$
has a local minimum at $r_Q = (1-Q^2)^{1/2}$ and a local maximum at
$r_q = (1-q^2)^{1/2}$, where $q,Q$ are implicitely defined by
$J = q(1-q^2) = Q(1-Q^2)$ and $q^2 < 1/3 < Q^2$. \cr

In the limiting case $J^2 = 4/27$, the effective potential
$V_J$ is stricly decreasing over $\R_+$ with an inflexion
point at $r = (2/3)^{1/2}$. Since $W$ is bounded, we must have
$E = E_-(J) = E_+(J) = 1/3$, hence (up to a global phase)
$W(x) = (1-q^2)^{1/2}\rme^{\rmi qx}$ with $q = (1/3)^{1/2}\sign(J)$.

Finally, if $W$ is a solution of \reff{e:sgl} with $J = 0$,
then (replacing $W(x)$ with $W(x)\,\rme^{\rmi \phi}$ for some
$\phi \in [0,2\pi])$ we can assume that $W(x) \in \R$ for all
$x \in \R$. Then $E = \frac12 r_x^2 + V_0(r)$ and we have the
same discussion as above with $E_-(0) = 0$ and $E_+(0) = 1/4$.
If $E = 0$, then $W \equiv 0$. If $E = 1/4$, then either
$W(x) \equiv \pm 1$ or $W(x) = \pm\tanh(x/\sqrt{2})$ (up to
a translation). If $0 < E < 1/4$, then $W$ is periodic with
{\it half-period}
\[
  T(0,E) \,=\, 2\int_0^{r_2(0,E)} \frac{\dd r}
  {\sqrt{2(E - V_0(r))}}~, \quad \hbox{where}\quad
  r_2(0,E) = (1-\sqrt{1-4E})^{1/2}~.
\]
Moreover, $W(x+T(0,E)) = -W(x)$, hence $|W|$ is periodic with
(minimal) period $T(0,E)$.
\QED

\medskip
As we shall see at the end of this section,
Proposition~\ref{th:glstat} implies the existence of a
six-parameter family of quasi-periodic solutions of the nonlinear
Schr\"odinger equation \reff{e:nls}. Before doing that, we study
in detail the properties of the period $T$ and the phase increment
$\Phi$ defined in \reff{e:Tdef}, \reff{e:Phidef}, because these
quantities play a crucial role in the stability analysis of the
solutions of \reff{e:sgl}, both for the Schr\"odinger and the
Ginzburg-Landau dynamics. We first give explicit formulas for $T$
and $\Phi$ (see also \cite{DGJ,EGW}) which are convenient for
analytical study and numerical approximation.

\begin{Lemma}\label{th:TPhi}
Assume that $(J,E) \in D$ and denote by $0 \le y_1 < y_2 < y_3$ the
roots of the cubic polynomial $P(y) = y^3 - 2y^2 + 4Ey - 2J^2$.
Then
\begin{equation}\label{e:Texp}
  T(J,E) \,=\, \sqrt{2} \int_{y_1}^{y_2} \frac{\dd y}{\sqrt{(y-y_1)
  (y-y_2)(y-y_3)}} \,=\, 2\sqrt{2} \int_0^{\pi/2}
  \frac{\dd \phi}{\sqrt{y_3 - s(\phi)}}~,
\end{equation}
where $s(\phi) = y_1 \cos^2(\phi) + y_2 \sin^2(\phi)$. Similarly,
if $J \neq 0$,
\begin{equation}\label{e:Phiexp}
  \Phi(J,E) \,=\, \sqrt{2} \int_{y_1}^{y_2} \frac{J}y\,
  \frac{\dd y}{\sqrt{(y-y_1)(y-y_2)(y-y_3)}} \,=\,
  2\sqrt{2} \int_0^{\pi/2} \frac{J}{s(\phi)}\,
  \frac{\dd \phi}{\sqrt{y_3 - s(\phi)}}~.
\end{equation}
\end{Lemma}

\figurewithtex 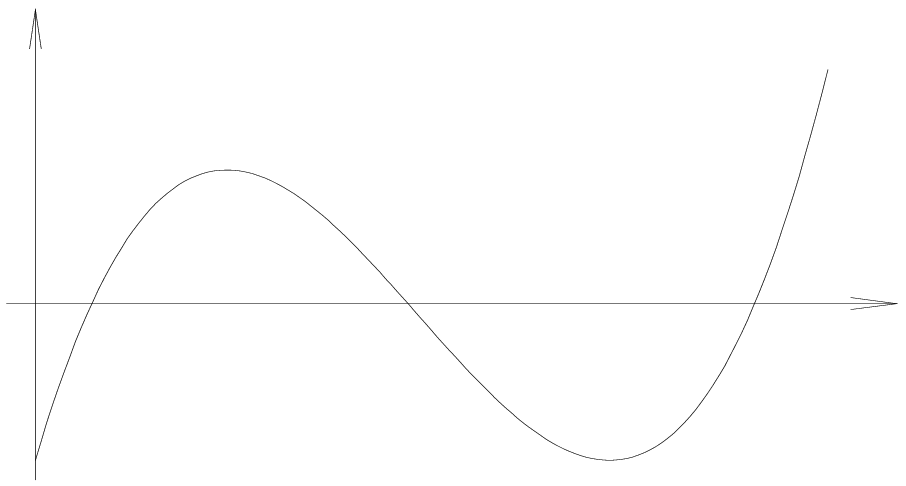 Fig3.tex 5.000 9.000
{\bf Fig.~3:} The three roots of the cubic polynomial $P(y) = y^3 - 2y^2
+4Ey - 2J^2$ for $(J,E) \in D$.\cr

\noindent{\bf Proof.} Observe that $P(r^2) = 4r^2 (E-V_J(r))$,
hence $P(y)$ has three nonnegative roots whenever $(J,E) \in D$
(see Fig.~2). Thus, using the change of variables $r = \sqrt{y}$
in \reff{e:Tdef}, we obtain the first expression in \reff{e:Texp}.
The last expression then follows by setting $y = s(\phi)$, so that
$\dd y = 2\sqrt{(y-y_1)(y_2-y)}\dd \phi$. Similarly, if $J \neq
0$, we obtain \reff{e:Phiexp} from \reff{e:Phidef}. \QED

\medskip
It follows in particular from \reff{e:Texp} that $T$ is a smooth
function of $(J,E) \in D$. In contrast $\Phi$ cannot be extended
to a continuous function over $D$, because (as is easily verified)
\[
  \lim_{J \to 0^\pm} \Phi(J,E) \,=\, \pm \pi~, \quad
  \hbox{for} \quad 0 < E < 1/4~.
\]
This suggests to introduce the renormalized phase $\Psi : D \to \R$
defined by
\begin{equation}\label{e:Psidef}
  \Psi(J,E) \,=\, \left\{
  \begin{array}{ccl}
  \Phi(J,E) - \pi\,\sign(J) & \hbox{if} & J \neq 0~,\\
  0 & \hbox{if} & J = 0~.
  \end{array}\right.
\end{equation}

\begin{Lemma}\label{th:Psi}
For any $(J,E) \in D$, one has
\begin{equation}\label{e:Psiexp}
  \Psi(J,E) \,=\, 2\sqrt{2} \,J \int_0^{\pi/2} \frac{\dd \phi}{\sqrt{y_3}
  \,\sqrt{y_3-s(\phi)}\,(\sqrt{y_3} + \sqrt{y_3-s(\phi)})}~.
\end{equation}
In particular, $\Psi : D \to \R$ is a smooth function.
\end{Lemma}

\noindent{\bf Proof.} The integral in \reff{e:Psiexp} is a smooth
function of $(J,E) \in D$, so it is sufficient to establish
\reff{e:Psiexp} for $J \neq 0$. In that case, we remark that
\[
  2\sqrt{2}\int_0^{\pi/2} \frac{J}{s(\phi)}\,\frac1{\sqrt{y_3}}
  \dd \phi \,=\, \frac{2\sqrt{2}}{\sqrt{y_3}} \Bigl(\frac{\pi}2
  \frac{J}{\sqrt{y_1y_2}}\Bigr) \,=\, \frac{\pi\sqrt{2}\,J}
  {\sqrt{y_1 y_2 y_3}} \,=\, \pi\,\sign(J)~,
\]
because $y_1 y_2 y_3 = 2 J^2$. Thus, using \reff{e:Phiexp} and
\reff{e:Psidef}, we obtain
\begin{eqnarray*}
  \Psi(J,E) &=& 2\sqrt{2} \int_0^{\pi/2}\frac{J}{s(\phi)}\,
  \Bigl(\frac1{\sqrt{y_3-s(\phi)}}-\frac1{\sqrt{y_3}}\Bigr)\dd \phi
  \\
  &=& 2\sqrt{2} J \int_0^{\pi/2} \frac{\dd \phi}{\sqrt{y_3}
  \,\sqrt{y_3-s(\phi)}\,(\sqrt{y_3} + \sqrt{y_3-s(\phi)})}~,
\end{eqnarray*}
which is the desired formula. \QED

\medskip
We next study the monotonicity properties of $T$ and $\Psi$. Since
$T(J,E)$ is an {\em even} and $\Psi(J,E)$ an {\em odd} function
of $J$, we may restrict ourselves to the half-domain $D_+ = \{(J,E)
\in D\,|\, J > 0\}$.

\begin{Proposition}\label{th:monotone}$ $\\[1mm]
i) $~\displaystyle \frac{\partial T}{\partial E}(J,E) > 0~$ for all
$(J,E) \in D$. \\[1mm]
ii) $~\displaystyle \frac{\partial \Psi}{\partial E}(J,E) \equiv
\frac{\partial \Phi}{\partial E}(J,E) =
-\frac{\partial T}{\partial J}(J,E) > 0~$ for all $(J,E) \in D_+$.
\end{Proposition}

\noindent{\bf Proof.} The monotonicity of the period $T$ with
respect to the energy $E$ has been established in \cite{DGJ}.
We give here a different argument, which is also a preparation
for the proof of Proposition~\ref{th:KAM} below. Let $(J,E) \in D$.
Since $y_1,y_2,y_3$ are solutions of the cubic equation
$y^3 - 2y^2 + 4Ey -2J^2 = 0$, we have $y_1 + y_2 + y_3 = 2$
and
\begin{equation}\label{e:yiEJ}
  \frac{\partial y_i}{\partial E} \,=\, -\frac{4y_i}{3y_i^2 - 4y_i
  + 4E}~, \quad
  \frac{\partial y_i}{\partial J} \,=\, \frac{4J}{3y_i^2 - 4y_i
  + 4E}~, \quad i = 1,2,3~.
\end{equation}
In particular,
\begin{equation}\label{e:yEsign}
  \frac{\partial y_1}{\partial E} \,<\, 0~, \quad
  \frac{\partial y_2}{\partial E} \,>\, 0~, \quad
  \frac{\partial y_1}{\partial E} + \frac{\partial y_2}{\partial E}
  \,=\, -\frac{\partial y_3}{\partial E} \,>\, 0~,
\end{equation}
because $P'(y_1) > 0$, $P'(y_2) < 0$, $P'(y_3) > 0$ (see Fig.~3).
Similarly, if $J > 0$,
\begin{equation}\label{e:yJsign}
  \frac{\partial y_1}{\partial J} \,>\, 0~, \quad
  \frac{\partial y_2}{\partial J} \,<\, 0~, \quad
  \frac{\partial y_1}{\partial J} + \frac{\partial y_2}{\partial J}
  \,=\, -\frac{\partial y_3}{\partial J} \,<\, 0~.
\end{equation}
On the other hand, differentiating \reff{e:Texp} with respect to
$E$ and $J$, we find
\begin{equation}\label{e:Tder}
  \frac{\partial T}{\partial E} \,=\, A_1 \frac{\partial y_1}{\partial E}
  + A_2 \frac{\partial y_2}{\partial E}~, \quad
  \frac{\partial T}{\partial J} \,=\, A_1 \frac{\partial y_1}{\partial J}
  + A_2 \frac{\partial y_2}{\partial J}~, \quad
\end{equation}
where
\[
  A_1 \,=\, \sqrt{2}\int_0^{\pi/2} \frac{1+\cos^2(\phi)}{(y_3-s(\phi))^{3/2}}
  \dd\phi~, \quad
  A_2 \,=\, \sqrt{2}\int_0^{\pi/2} \frac{1+\sin^2(\phi)}{(y_3-s(\phi))^{3/2}}
  \dd\phi~. \quad
\]
The crucial observation is:
\begin{equation}\label{e:A1A2}
  A_2 - A_1 \,=\, \sqrt{2}\int_0^{\pi/2}\frac{\sin^2(\phi)-\cos^2(\phi)}
  {(y_3-s(\phi))^{3/2}}\dd\phi \,>\, 0~.
\end{equation}
Indeed, this inequality follows from Lemma~\ref{th:fg} below,
with $I = [0,\pi/2]$, $\dd\mu = (2/\pi)\dd \phi$, $f(\phi) =
\sin^2(\phi)-\cos^2(\phi)$ and $g(\phi) = (y_3-s(\phi))^{-3/2}$
(remark that $f,g$ are strictly increasing, and that $\int_0^{\pi/2}
f(\phi)\dd\phi = 0$.) Thus $A_2 > A_1 > 0$. We conclude that
\begin{equation}\label{e:TEpos}
  \frac{\partial T}{\partial E} \,=\, (A_2-A_1)\frac{\partial y_2}
  {\partial E} + A_1 \Bigl(\frac{\partial y_1}{\partial E} +
  \frac{\partial y_2}{\partial E}\Bigr) \,>\, 0~.
\end{equation}
Similarly, if $J > 0$,
\[
  \frac{\partial T}{\partial J} \,=\, (A_2-A_1)\frac{\partial y_2}
  {\partial J} + A_1 \Bigl(\frac{\partial y_1}{\partial J} +
  \frac{\partial y_2}{\partial J}\Bigr) \,<\, 0~.
\]

Finally, the relation $\partial T/\partial J = -\partial \Phi/
\partial E$ is a consequence of the following standard observation
\cite{BR1}: if $\AA : D \to \R$ is the {\em action functional}
defined by
\begin{equation}\label{e:AAdef}
  \AA(J,E) \,=\, 2\int_{r_1(J,E)}^{r_2(J,E)} \sqrt{2(E - V_J(r))}
  \dd r~,
\end{equation}
then a direct calculation shows that $T = \partial \AA/\partial E$
and $\Phi = -\partial\AA/\partial J$, and the result follows. \QED

\medskip
The following elementary result will be used several times in this
paper:

\begin{Lemma}\label{th:fg}
Let $\mu$ be a (Borel) probability measure on some interval
$I \subset \R$, and let $f,g : I \to \R$ be bounded and measurable
functions. Then
\[
  \int_I f(x)g(x)\dd\mu - \Bigl(\int_I f(x)\dd\mu\Bigr)
  \Bigl(\int_I g(x)\dd\mu\Bigr) \,=\, \frac12
  \int_{I\times I} (f(x)-f(y))(g(x)-g(y))\dd\mu_x \dd\mu_y~.
\]
In particular, if both $f$ and $g$ are strictly increasing or
strictly decreasing, and if the support of $\mu$ is not reduced
to a single point, then
\[
  \int_I f(x)g(x)\dd\mu  \,>\, \Bigl(\int_I f(x)\dd\mu\Bigr)
  \Bigl(\int_I g(x)\dd\mu\Bigr)~.
\]
\end{Lemma}

We now establish an important non-degeneracy property of
system \reff{e:sgl}. For any $(J,E) \in D$, let
\begin{equation}\label{e:Deltadef}
  \Delta(J,E) \,=\, \Det\left(\begin{array}{cc}
  \frac{\partial T}{\partial E} &
  \frac{\partial \Psi}{\partial E} \\[1mm]
  \frac{\partial T}{\partial J} &
  \frac{\partial \Psi}{\partial J}
  \end{array}\right)(J,E)~.
\end{equation}

\begin{Proposition}\label{th:KAM}
For all $(J,E) \in D$ we have $\Delta(J,E) > 0$.
\end{Proposition}

\noindent{\bf Remarks}\\
{\bf 1.} Proposition~\ref{th:KAM} shows that the integrable
Hamiltonian system \reff{e:sgl} is non-degenerate in the
sense of KAM theory. Indeed, by Liouville's theorem, we can
express (at least locally) the energy $E$ of the system in
terms of the {\em action variables} $J$ and $\AA$, where
$\AA$ is defined in \reff{e:AAdef}. If we denote $E = H(J,\AA)$,
a direct calculation \cite{BR1} shows that
\[
  \Delta \,=\, -T^4 \,\Det\left(\begin{array}{cc}
  \frac{\partial^2 H}{\partial \AA^2} &
  \frac{\partial^2 H}{\partial \AA\,\partial J} \\[1mm]
  \frac{\partial^2 H}{\partial J\,\partial \AA} &
  \frac{\partial^2 H}{\partial J^2}
  \end{array}\right)~,
\]
hence the KAM determinant (the Hessian of $H$) is always negative.\\[1mm]
{\bf 2.} By Proposition~\ref{th:monotone}, it is clear that
$\Delta > 0$ whenever $\partial \Psi/\partial J > 0$. Unfortunately,
the latter inequality is not true for all $(J,E) \in D$.

\medskip\noindent{\bf Proof.}
To evaluate $\Delta$ we have to compute the derivatives of the
renormalized phase $\Psi$ with respect to $E$ and $J$. Using
\reff{e:Psiexp}, we obtain after straightforward calculations:
\begin{equation}\label{e:Psider}
  \frac{\partial \Psi}{\partial E} \,=\, B_1 J\frac{\partial y_1}{\partial E}
  + B_2 J\frac{\partial y_2}{\partial E}~, \quad
  \frac{\partial \Psi}{\partial J} \,=\, B_1 J\frac{\partial y_1}{\partial J}
  + B_2 J\frac{\partial y_2}{\partial J} + B_3~, \quad
\end{equation}
where
\begin{eqnarray}\nonumber
  {B_1 \atop B_2} &=&
  \sqrt{2} \int_0^{\pi/2} \biggl\{\frac{1}{y_3^{3/2}(y_3-s(\phi))^{1/2}}
  \,\frac{2\sqrt{y_3}+\sqrt{y_3-s(\phi)}}{(\sqrt{y_3}+\sqrt{y_3-s(\phi)})^2}
  \,+ \\ \label{e:Bidef}
  && \qquad \qquad \frac{1}{y_3^{1/2}(y_3-s(\phi))^{3/2}}
  \,\frac{\sqrt{y_3}+2\sqrt{y_3-s(\phi)}}{(\sqrt{y_3}+\sqrt{y_3-s(\phi)})^2}
  \,\Bigl(1+{\cos^2(\phi) \atop \sin^2(\phi)}\Bigr)\biggr\}\dd\phi~,
  \\[1mm]\nonumber
  B_3 &=& 2\sqrt{2} \int_0^{\pi/2} \frac{1}{\sqrt{y_3}
  \,\sqrt{y_3-s(\phi)}\,(\sqrt{y_3} + \sqrt{y_3-s(\phi)})}\dd\phi
  \,=\, \frac{\Psi}{J}~.
\end{eqnarray}
(Here the formula for $B_1$ should be read with $\cos^2(\phi)$ in
the right-hand side, and $B_2$ with $\sin^2(\phi)$; both expressions
are otherwise identical.) Using \reff{e:Tder}, \reff{e:Psider}
together with the identities
\[
  J \frac{\partial y_i}{\partial E} + y_i \frac{\partial y_i}
  {\partial J} \,=\, 0~, \quad i = 1,2,3~,
\]
which follow from \reff{e:yiEJ}, we find
\begin{eqnarray}\nonumber
 \Delta(J,E) &=& \Bigl(A_1 \frac{\partial y_1}{\partial E}
  + A_2 \frac{\partial y_2}{\partial E}\Bigr)
  \Bigl(B_1 J\frac{\partial y_1}{\partial J}
  + B_2 J\frac{\partial y_2}{\partial J} + B_3\Bigr) \\\nonumber
  &&-\,\Bigl(A_1 \frac{\partial y_1}{\partial J}
  + A_2 \frac{\partial y_2}{\partial J}\Bigr)
  \Bigl(B_1 J\frac{\partial y_1}{\partial E}
  + B_2 J\frac{\partial y_2}{\partial E}\Bigr)\\\nonumber
 &=& J(A_1 B_2 - A_2 B_1)\Bigl(\frac{\partial y_1}{\partial E}
   \frac{\partial y_2}{\partial J} - \frac{\partial y_2}{\partial E}
   \frac{\partial y_1}{\partial J}\Bigr) + B_3
   \frac{\partial T}{\partial E}\\\label{e:Deltaexp}
 &=& (A_1 B_2 - A_2 B_1)(y_2-y_1)\frac{\partial y_1}{\partial J}
   \frac{\partial y_2}{\partial J} + B_3
   \frac{\partial T}{\partial E}~.
\end{eqnarray}
Since $\partial T/\partial E > 0$, $B_3 > 0$, $y_2-y_1 > 0$ and
$(\partial y_1/\partial J)(\partial y_2/\partial J) < 0$, it is
sufficient to prove that $\Delta_1 \eqdef A_2 B_1 - A_1 B_2 > 0$.

To achieve this goal, we set $\sigma(\phi) =
(1-y_3^{-1}s(\phi))^{1/2}$ and we observe that $A_i =
\sqrt{2}\,y_3^{-3/2}\bA_i$, $B_i = \sqrt{2} \,y_3^{-5/2}\bB_i$,
for $i = 1,2$, where
\begin{eqnarray}\label{e:bAiBidef}
  \bA_1 &=& \int_0^{\pi/2} \frac{1+\cos^2(\phi)}{\sigma(\phi)^3}\dd\phi~,
  \quad \bA_2 ~=~ \int_0^{\pi/2} \frac{1+\sin^2(\phi)}{\sigma(\phi)^3}
  \dd\phi~,\\[1mm]\nonumber
  {\bB_1 \atop \bB_2} &=& \int_0^{\pi/2}\biggl\{\frac{2+\sigma(\phi)}
  {\sigma(\phi)(1+\sigma(\phi))^2} + \frac{1+2\sigma(\phi)}
  {\sigma(\phi)^3(1+\sigma(\phi))^2}\Bigl(1+{\cos^2(\phi) \atop
  \sin^2(\phi)}\Bigr)\biggr\}\dd\phi~.
\end{eqnarray}
Thus $\Delta_1 = 2y_3^{-4}\Delta_2$, where $\Delta_2 = \bA_2\bB_1
-\bA_1\bB_2$. To prove that $\Delta_2 > 0$, we remark that
\[
  \frac{1+2\sigma}{\sigma^3(1+\sigma)^2} \,=\, \frac{1}{\sigma^3}
  - \frac{1}{\sigma(1+\sigma)^2}~, \quad \hbox{and}\quad
  \frac{2+\sigma}{\sigma(1+\sigma)^2} - \frac{1+\cos^2(\phi)}
  {\sigma(1+\sigma)^2} \,=\, \frac{\sigma+\sin^2(\phi)}
  {\sigma(1+\sigma)^2}~.
\]
It follows that $\bB_1 = \bA_1 + \tilde \bB_1$ and $\bB_2 = \bA_2 +
\tilde \bB_2$, where
\begin{equation}\label{e:tildeBidef}
  \tilde \bB_1 \,=\, \int_0^{\pi/2}\frac{\sigma(\phi)+\sin^2(\phi)}
  {\sigma(\phi)(1+\sigma(\phi))^2}\dd\phi~, \quad
  \tilde \bB_2\,=\, \int_0^{\pi/2}\frac{\sigma(\phi)+\cos^2(\phi)}
  {\sigma(\phi)(1+\sigma(\phi))^2}\dd\phi~.
\end{equation}
Now, we know from \reff{e:A1A2} that $\bA_2 > \bA_1 > 0$, and the
same argument (using Lemma~\ref{th:fg}) shows that $\tilde \bB_1 >
\tilde \bB_2 > 0$. Thus
\[
  \Delta_2 \,=\, \bA_2\bB_1 - \bA_1\bB_2 \,=\, \bA_2 \tilde\bB_1
  -\bA_1 \tilde\bB_2 \,>\, 0~,
\]
and the proof is complete. \QED

\figurewithtex 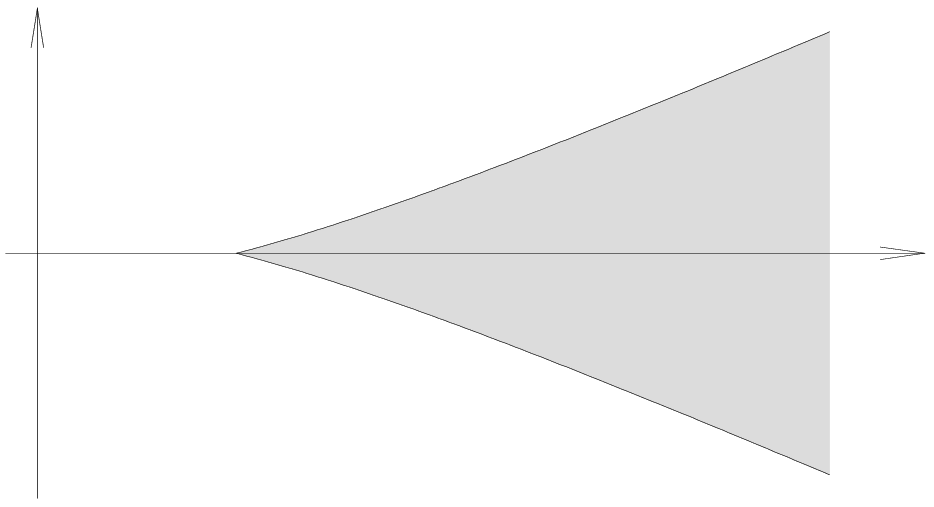 Fig4.tex 5.000 10.000
{\bf Fig.~4:} The region $\hat D \subset \R^2$ consisting of all
possible values of the pair $(T,\Psi)$.\cr

\medskip
Propositions~\ref{th:monotone} and \ref{th:KAM} imply that the
quasi-periodic solutions of \reff{e:sgl} can be parametrized
by the period $T$ and the renormalized phase $\Psi$, instead
of the angular momentum $J$ and the energy $E$. Indeed, let
\[
  \hat D \,=\, \Big\{(T,\Psi) \in \R^2 \,\Big|\, T > \pi\,,~
  |\Psi| < \hat \Psi(T)\Big\}~, \quad \hbox{where} \quad
  \hat \Psi(T) \,=\, \Bigl(\frac{T^2+2\pi^2}{3}\Bigr)^{1/2} - \pi~.
\]

\begin{Proposition}\label{th:Tpsirange}
The map $(T,\Psi) : D \to \hat D$ is a smooth diffeomorphism.
\end{Proposition}

\noindent{\bf Proof.} Let $J = Q(1-Q^2)$, where $1/3 < Q^2 \le 1$.
Using \reff{e:Texp} it is straightforward to verify that
\begin{equation}\label{e:Tbound}
  \lim_{E \to E_-(J)} T(J,E) \,=\, \frac{\pi\sqrt{2}}{\sqrt{3Q^2 - 1}}
    \,\equiv\, \frac{\pi}{(1-3E_-(J))^{1/4}}~, \quad
  \lim_{E \to E_+(J)} T(J,E) \,=\, +\infty~.
\end{equation}
Similarly, if $1/\sqrt{3} < Q < 1$ (so that $J > 0$), we deduce from
\reff{e:Phiexp} that
\begin{equation}\label{e:Phibound}
  \lim_{E \to E_-(J)} \Phi(J,E) \,=\, \frac{\pi\sqrt{2}\,Q}
  {\sqrt{3Q^2 - 1}}~, \quad
  \lim_{E \to E_+(J)} \Phi(J,E) \,=\, +\infty~.
\end{equation}
Since $\partial T/\partial E > 0$ by Proposition~\ref{th:monotone},
it follows from \reff{e:Tbound} that the range of the map $T : D \to
\R$ is exactly the interval $(\pi,+\infty)$. Fix $T_0 > \pi$
and let $\Sigma = \{(J,E) \in D\,|\, T(J,E) = T_0\}$. By the
Implicit Function Theorem, $\Sigma$ is a smooth curve in $D$ which
can be represented as a graph over the $J$-axis. Moreover, we know
that $\Sigma$ connects the boundary points $(-J_0,E_-(J_0))$ and
$(J_0,E_-(J_0))$, where $J_0 = Q_0(1-Q_0^2)$ and $Q_0 \in
(1/\sqrt{3},1)$ is determined by the relation
\[
  \frac{\pi\sqrt{2}}{\sqrt{3Q_0^2-1}} \,=\, T_0~,
\]
(see \reff{e:Tbound}). Now, Proposition~\ref{th:KAM} implies that
the restriction of $\Psi$ to the curve $\Sigma$ is a strictly
increasing function of $J$, because
\[
  \frac{\rmd }{\rmd \,J} \Psi|_\Sigma \,=\,
  \Bigl(\frac{\partial T}{\partial E}\Bigr)^{-1}
  \Bigl(\frac{\partial T}{\partial E}\frac{\partial \Psi}{\partial J}
  - \frac{\partial T}{\partial J}\frac{\partial \Psi}{\partial E}
  \Bigr) \,>\, 0~.
\]
Thus $\Psi$ varies from $-\Psi_0$ to $\Psi_0$ on the curve $\Sigma$,
where by \reff{e:Phibound}
\[
  \Psi_0 \,=\, \frac{\pi\sqrt{2}\,Q_0}{\sqrt{3Q_0^2-1}} - \pi
  \,=\, \Bigl(\frac{T_0^2+2\pi^2}{3}\Bigr)^{1/2} - \pi \,\equiv\,
  \hat\Psi(T_0)~.
\]
This proves that $(T,\Psi) : D \to \hat D$ is onto, and the
monotonicity properties established in Proposition~\ref{th:monotone}
imply that $(T,\Psi)$ is also one-to-one. \QED

\medskip
To conclude this section, we briefly verify that
Proposition~\ref{th:glstat} implies the existence of periodic waves
of the NLS equation of the form \reff{e:qpdef}. Fix $(J,E) \in D$,
and let $W : \R \to \C$ be a bounded solution of \reff{e:sgl}
satisfying \reff{e:JEdef}. (The proof of Proposition~\ref{th:glstat}
shows that this solution is unique up to a translation and a phase
factor. If needed, we can specify a particular solution by imposing
for instance $W(0) = r_2$, $W'(0) = \rmi J/r_2$, where $r_2 > 0$ is as
in Fig.~2.) We now set
\begin{equation}\label{e:Prep}
  W(x) \,=\, \rme^{\rmi \ell x}\,P(kx)~, \quad x \in \R~,
\end{equation}
where
\begin{equation}\label{e:kldef}
  k \,=\, \frac{\pi}{T(J,E)}~, \quad \hbox{and}\quad
  \ell \,=\, \frac{\Psi(J,E)}{T(J,E)}~.
\end{equation}
As $x \mapsto |W(x)|$ is periodic with minimal period $T(J,E)$, it is
clear that $y \mapsto |P(y)|$ is periodic with minimal period $\pi$.
Moreover, since $W(x+T) = \rme^{\rmi \Phi}W(x)$ by definition
of $\Phi(J,E)$, we also have $P(y+\pi) = -P(y)$ for all $y \in \R$,
hence $P$ is $2\pi$-periodic. Thus $U(x,t) = \rme^{-\rmi t}W(x) =
\rme^{\rmi (\ell x - t)}P(kx)$ is a quasi-periodic solution
of \reff{e:nls} of the form \reff{e:qpdef}, with $\omega = 1$
and $c = 0$.

\begin{Remark}
Using the continuous symmetries of the NLS equation, we can produce
for each pair $(J,E) \in D$ a four-parameter family of
periodic waves:
\[
  U_{c,\lambda,\phi,\xi}(x,t) \,=\, \lambda \,\rme^{\rmi (p_{c,\lambda}x
  -\omega_{c,\lambda}t - \phi)}P(k\lambda(x-ct)+\xi)~,
\]
where $c \in \R$, $\lambda > 0$, $\phi,\xi \in [0,2\pi]$, and
$p_{c,\lambda} = \lambda\ell + c/2$, $\omega_{c,\lambda} =
\lambda^2 + c\lambda\ell +c^2/4$. Taking into account the
parameters $J,E$, we obtain altogether a six-parameter family of
periodic waves of \reff{e:nls}.
\end{Remark}

The representation \reff{e:Prep} is well-adapted to understand the
connection between the bounded solutions $W$ of \reff{e:sgl} and the
periodic waves of the NLS equation, especially in the case of the cnoidal
waves for which $J = \ell = 0$. However it is more convenient for our
purposes to write the solution $W$ of \reff{e:sgl} in the alternative
form
\begin{equation}\label{e:Qrep}
  W(x) \,=\, \rme^{\rmi (\ell + k)x}\,Q^+(2kx) \,=\, \rme^{\rmi
  (\ell - k)x}\,Q^-(2kx)~, \quad x \in \R~,
\end{equation}
where $Q^\pm(z) = \rme^{\mp \rmi z/2}P(z/2)$. By construction, $Q^\pm$
and $|Q^\pm|$ are now periodic functions with the {\it same} minimal
period $2\pi$. This property facilitates the description of the
special class of perturbations that we use for the stability analysis,
see the statement of Theorem~\ref{th:orbit} and the remarks
thereafter. The representation \reff{e:Qrep} is also very natural for
solutions which are close to plane waves: in such a situation, 
either $Q^+$ or $Q^-$ is close to a constant, depending on the sign 
of $J$. This follows from the fact that $\ell+k = \Phi/T$ if $J > 0$ and 
$\ell-k = \Phi/T$ if $J < 0$.

\section{Orbital stability}
\label{s:os}

Our aim in this section is to show that the periodic waves of the
defocusing nonlinear Schr\"odinger equation \reff{e:nls} with
$\gamma = -1$ are stable within the class of solutions which have
the same period $T$ and the same Floquet multiplier $\rme^{\rmi
\Phi}$. Given $(J,E) \in D$, where $D \subset \R^2$ is the
parameter domain \reff{e:Ddef}, we consider the periodic wave
$U_{J,E}(x,t) = \rme^{-\rmi t} W_{J,E}(x)$, where $W_{J,E}$ is the
unique solution of \reff{e:sgl} with initial data $W_{J,E}(0) =
r_2$, $W_{J,E}'(0) = \rmi J/r_2$, and $r_2 > 0$ is as in Fig.~2.
In particular, $W_{J,E}$ satisfies \reff{e:JEdef}.  Let $p =
k+\ell$, where $k \equiv k_{J,E}$ and $\ell \equiv \ell_{J,E}$ are
defined by \reff{e:kldef}. As in \reff{e:Qrep}, we set $W_{J,E}(x)
= \rme^{\rmi px} Q_{J,E}(2kx)$ and we recall that $Q_{J,E}(z)$ is
a $2\pi$-periodic function of $z$.

To study the stability of the periodic wave $U_{J,E}(x,t)$, we
consider solutions of \reff{e:nls} of the form
\begin{equation}\label{e:Qgen}
  U(x,t) \,=\, \rme^{\rmi (px-t)} Q(2kx,t)~,
\end{equation}
where $Q(z,t)$ satisfies the evolution equation
\begin{equation}\label{e:Qeq}
  \rmi Q_t + 4\rmi pk Q_z + 4k^2 Q_{zz} + (1-p^2)Q - |Q|^2Q \,=\, 0~.
\end{equation}
By construction, $Q_{J,E}(z)$ is now a stationary solution of
\reff{e:Qeq} and our goal is to show that this equilibrium is
stable with respect to $2\pi$-periodic perturbations. We thus
introduce the function space
\[
  X \,=\, H^1_\per([0,2\pi],\C) \,=\, \Bigl\{u \in H^1_\loc(\R,\C)
  \,\Big|\, u(z) = u(z+2\pi) \hbox{ for all } z \in \R\Bigr\}~,
\]
which is viewed as a {\it real} Hilbert space equipped with the
scalar product
\[
  (u,v)_X \,=\, \Re \int_0^{2\pi} (u(z) \overline{v}(z) +
  u_z(z) \overline{v}_z(z)) \dd z~, \quad u,v \in X~.
\]
As usual, the dual space $X^*$ will be identified with
$H^{-1}_\per([0,2\pi],\C)$ through the pairing
\[
  \langle u,v \rangle \,=\, \Re \int_0^{2\pi} u(z) \overline{v}(z)
  \dd z~, \quad u \in X^*~, \quad v \in X~.
\]

It is well-known that the Cauchy problem for \reff{e:Qeq} is
globally well-posed in the space $X$ (see \cite{CW,GV1,GV2,K2}).
Moreover, the evolution defined by \reff{e:Qeq} on $X$ is
invariant under the action of a two-parameter group of isometries:
the space translations and the phase rotations. The symmetry group
$G$ is thus the two-dimensional torus $\T^2 = (\R / 2\pi\Z)^2$
which acts on $X$ through the unitary representation $\RR$ defined
by
\[
  (\RR_{(\phi,\xi)} u)(z) \,=\, \rme^{-\rmi \phi} u(z+\xi)~,
  \quad u \in X~, \quad (\phi,\xi) \in G~.
\]
In view of these symmetries, it is natural to introduce the
semi-distance $\rho$ on $X$ defined by
\begin{equation}\label{e:rhodef}
  \rho(u,v) \,=\, \inf_{(\phi,\xi) \in G} \|u - \RR_{(\phi,\xi)}v\|_X~,
  \quad u,v \in X~.
\end{equation}
The main result of this section is the following reformulation of
Theorem~\ref{th:orbit}.

\begin{Proposition}\label{th:Qorbit}
Given $(J,E) \in D$, there exist $C_0 > 0$ and $\epsilon_0 > 0$
such that, if $Q_0 \in X$ satisfies $\rho(Q_0,Q_{J,E}) \le \epsilon$
for some $\epsilon \le \epsilon_0$, then the solution $Q(z,t)$
of \reff{e:Qeq} with initial data $Q_0$ satisfies $\rho(Q(\cdot,t),
Q_{J,E}) \le C_0 \epsilon$ for all $t \in \R$.
\end{Proposition}

To prove Proposition~\ref{th:Qorbit}, we follow the general approach
of Grillakis, Shatah, and Strauss \cite{GSS1,GSS2}. We first observe
that Eq.\reff{e:Qeq} inherits from the original NLS equation
several conserved quantities:
\begin{eqnarray}\nonumber
  && N(Q) \,=\, \frac12 \int_0^{2\pi} |Q(z)|^2 \dd z~, \\[2mm]
  \label{e:NMEdef}
  && M(Q) \,=\, \frac{\rmi}{2} \int_0^{2\pi} \overline Q(z)
     Q_z(z) \dd z~, \\[2mm] \nonumber
  && \EE(Q) \,=\, \int_0^{2\pi} \Bigl( 2k^2 |Q_z(z)|^2
     + \frac14|Q(z)|^4\Bigr)\dd z~,
\end{eqnarray}
which will be referred to as the {\em charge}, the {\em
momentum}, and the {\em energy}, respectively. The charge $N$
is conserved due to the phase invariance of \reff{e:Qeq}, the
momentum $M$ due to the translation invariance, and the energy
$\EE$ because \reff{e:Qeq} is autonomous. Clearly $N$, $M$, and $\EE$
are smooth real functions on $X$. Their first order derivatives
are therefore smooth maps from $X$ into $X^*$:
\begin{equation}\label{e:firstder}
  N'(Q) \,=\, Q~, \quad M'(Q) \,=\, \rmi Q_z~, \quad
  \EE'(Q) \,=\, -4k^2 Q_{zz} + |Q|^2 Q~.
\end{equation}
Similarly, the second order derivatives are smooth maps from
$X$ into $\mathcal{L}(X,X^*)$, the space of all bounded linear
operators from $X$ into $X^*$:
\[
  N''(Q) = \1~, \quad M''(Q) \,=\, \rmi \partial_z~, \quad
  \EE''(Q) \,=\, -4k^2\partial_{zz} + |Q|^2 + 2 Q\otimes Q~,
\]
where $\langle (Q\otimes Q)u,v\rangle = \int_0^{2\pi}
\Re(Q\overline{u})\Re(Q\overline{v})\dd z$ for all $u,v \in X$.

By construction, the periodic wave profile $Q_{J,E}$ is a critical
point of the modified energy
\begin{equation}\label{e:Eabdef}
  \EE_{J,E}(Q) \,=\, \EE(Q) - (1-p^2)N(Q) - 4p k M(Q)~,
\end{equation}
namely $\EE_{J,E}'(Q_{J,E}) = 0$. To determine the nature of
this critical point, we consider the second variation
\[
  H_{J,E} \,=\, \EE_{J,E}''(Q_{J,E}) \,=\, -4k^2 \partial_{zz}
  -4\rmi pk\partial_z - (1{-}p^2) + |Q_{J,E}|^2 + 2 Q_{J,E}
  \otimes Q_{J,E}~.
\]
Since $X$ is a real Hilbert space, it is natural to decompose
its elements (which are complex functions) into real and
imaginary parts, in which case we obtain the matrix operator
\[
  \hat H_{J,E} \,=\, \pmatrix{-4k^2\partial_{zz} -
  (1{-}p^2)+ 3 R_{J,E}^2 + I_{J,E}^2 & 4pk \partial_z
  + 2R_{J,E} I_{J,E}\cr - 4pk\partial_z + 2R_{J,E}I_{J,E}&
  -4k^2\partial_{zz} - (1{-}p^2) + R_{J,E}^2 + 3I_{J,E}^2}~,
\]
where $Q_{J,E} = R_{J,E} + \rmi I_{J,E}$. As is easily verified,
$H_{J,E}$ is a self-adjoint operator in
$X_0:=L^2_\per([0,2\pi],\C)$ with compact resolvent, and $H_{J,E}$
is bounded from below. According to \cite{GSS1, GSS2}, a crucial
information is the number of negative eigenvalues of this
operator.

\begin{Proposition}\label{th:Hneg}
For any $(J,E) \in D$, the operator $H_{J,E}$ acting on $X_0$ has
a simple negative eigenvalue, a double eigenvalue at zero, and the
rest of the spectrum is strictly positive.
\end{Proposition}

\noindent{\bf Proof.} Since Eq.\reff{e:Qeq} is invariant under the
action of $G$, we know that $\RR_{(\phi,\xi)} Q_{J,E}$ is a
stationary solution for all $(\phi,\xi) \in G$. Differentiating
with respect to $\phi$ and $\xi$ at $(\phi,\xi) = (0,0)$, we
obtain $H_{J,E} Q_{J,E}' = H_{J,E}(\rmi Q_{J,E}) = 0$. It is clear
that the functions $Q_{J,E}'$ and $\rmi Q_{J,E}$ are linearly
independent in the domain $X_2 := H^2_\per([0,2\pi],\C)$ of
$H_{J,E}$ (otherwise $W_{J,E}(x) = \rme^{\rmi px} Q_{J,E}(2kx) $
would be a plane wave, contradicting the assumption $(J,E) \in
D$). Thus we see that zero is an eigenvalue of $H_{J,E}$ of
multiplicity at least two.

Following \cite{DGJ}, we next show that the multiplicity of zero
as an eigenvalue of $H_{J,E}$ in $X_0$ is always {\em exactly
two}. The idea is to produce two other solutions of the
differential equation $H_{J,E} Q = 0$ by differentiating the
profile $Q_{J,E}$ with respect to $J$ and $E$, and to verify that
none of these solutions belongs to the domain $X_2$. Since the
kernel $K = \{Q \in C^2(\R,\C)\,|\, H_{J,E}Q = 0\}$ is a
four-dimensional (real) vector space, it will follow that $K \cap
X_2 = \linspan\{Q_{J,E}',\rmi Q_{J,E}\}$, which is the desired
result.

However, the profile $Q_{J,E}$ is a stationary solution of the
equation \reff{e:Qeq} which has coefficients depending upon $J$
and $E$, so that we cannot find the two solutions of $H_{J,E} Q =
0$ just by differentiating $Q_{J,E}$ with respect to $J$ and $E$.
Instead, we start with the quasi-periodic solution $W_{J,E}(x) =
\rme^{\rmi px} Q_{J,E}(2kx) $ of \reff{e:sgl}. In view of
Proposition~\ref{th:Tpsirange}, we can use $T$ and $\Psi$ (instead
of $J$ and $E$) to parametrize the family of quasi-periodic
solutions of \reff{e:sgl}. For our present purposes, the most
convenient set of parameters will be $(k,p)$, where $p = k+\ell$
and $k,\ell$ are defined in \reff{e:kldef}. Since the equation
\reff{e:sgl} does not depend upon $k$ and $p$, the derivatives
${\partial W_{J,E}}/{\partial k}$ and ${\partial
W_{J,E}}/{\partial p}$ are solutions of the linear equation
\[
  \omega'' + \omega -2|W_{J,E}|^2 \omega - W_{J,E}^2 \,
  \overline{\omega} \,=\, 0~.
\]
It is then straightforward to check that the functions $R_1(z)$
and $R_2(z)$ defined for $z=2kx$ through
\begin{eqnarray}\label{e:omidef}
  R_1(2kx) & = & \rme^{-\rmi px} \frac{\partial W_{J,E}}{\partial k}(x)
  \; =\; \frac{\partial Q _{J,E}}{\partial k}(2kx) + 2x Q_{J,E}'(2kx)~,
  \\[1ex]
  R_2(2kx) & = & \rme^{-\rmi px} \frac{\partial W_{J,E}}{\partial p}(x)
  \; =\; \frac{\partial Q _{J,E}}{\partial p}(2kx) + \rmi x Q_{J,E}(2kx)~,
  \nonumber
\end{eqnarray}
satisfy  $H_{J,E}R_1 = H_{J,E}R_2 = 0$, namely $R_1, R_2 \in K$.
Of course $R_1,R_2$ depend on the point $(J,E) \in D$ where the
derivative is taken, but we omit this dependence for notational
simplicity.

Now, let $R(z) = A_1 R_1(z) + A_2 R_2(z)$ for some $A_1, A_2 \in
\R$. Using the definitions \reff{e:omidef} and the fact that
$Q_{J,E}$ is $2\pi$-periodic we find that $R$ satisfies the
periodicity conditions $R(2\pi) = R(0)$ and $R'(2\pi) = R'(0)$ if
and only if
\[
  2A_1 Q_{J,E}'(0) + \rmi A_2 Q_{J,E}(0) \,=\,0~, \quad \hbox{and}
  \quad 2A_1 Q_{J,E}''(0) + \rmi A_2 Q_{J,E}'(0) \,=\,0~.
\]
This linear system has a nontrivial solution $(A_1,A_2) \in \C^2$
if and only if
\begin{equation}\label{e:Qcond}
  Q_{J,E}'(0)^2 \,=\, Q_{J,E}(0) \,Q_{J,E}''(0)~.
\end{equation}
Using again the quasi-periodic solution $W_{J,E}(x)$ which
satisfies the initial conditions $W_{J,E}(0) = r_2$, $W_{J,E}'(0)
= \rmi J/r_2$ and $W_{J,E}''(0) = r_2^3 -r_2$, it is not difficult
to verify that \reff{e:Qcond} is equivalent to $J^2 = r_2^4 -
r_2^6$. But we know that $J^2 = r_Q^4 - r_Q^6 = r_q^4 - r_q^6$
with $r_Q < r_2 < r_q$ (see Fig.~2), hence $J^2 < r_2^4 - r_2^6$
so that \reff{e:Qcond} never holds. We conclude that $R = A_1 R_1
+ A_2 R_2$ belongs to $X_2$ only if $A_1 = A_2 = 0$. This shows
that the four functions $Q_{J,E}', \rmi Q_{J,E}, R_1, R_2$ are
linearly independent (over $\R$) and therefore form a basis of the
kernel $K$. Moreover, $K \cap X_2 = \linspan\{Q_{J,E}',\rmi
Q_{J,E}\}$ as expected.

It is now easy to conclude the proof of Proposition~\ref{th:Hneg}.
Indeed, since the eigenvalues of $H_{J,E}$ depend continuously on
$(J,E) \in D$, the fact that the zero eigenvalue has constant
multiplicity implies that the number of negative eigenvalues of
$H_{J,E}$ remains unchanged when $(J,E)$ varies over $D$.  It is
therefore sufficient to verify that the conclusion of
Proposition~\ref{th:Hneg} holds for {\em one} value of $(J,E) \in
D$. But as $(J,E) \to 0$ the wave profile $Q_{J,E}$ converges
uniformly to zero, hence $H_{J,E}$ converges to a constant
coefficient operator whose spectrum is easily determined by
Fourier analysis. Using a perturbation argument, it is then
straightforward to localize the eigenvalues of $H_{J,E}$ for
$(J,E)$ close to the origin. This calculation is performed in
\cite{GH1}, using an appropriate parametrization of the small
amplitude periodic waves. In the particular case $J = 0$, we
obtain the result below, which concludes the proof of
Proposition~\ref{th:Hneg}. \QED

\begin{Lemma}[{\cite[Proposition A.1]{GH1}}]\label{th:smallwaves}
If $J = 0$ and $E > 0$ is sufficiently small, the operator
$H_{J,E}$ acting on $X_0$ has exactly three eigenvalues
$\{\lambda_1,\lambda_2,0\}$ in a neighborhood of the origin, and
the rest of its spectrum is positive and bounded away from zero.
Moreover, the eigenvalues $\lambda_1,\lambda_2$ are simple and
satisfy \footnote{The expansions of $\lambda_1, \lambda_2$ do not 
appear explicitly in \cite{GH1}, but they follow easily from 
the proof of Proposition~A.1 there if one observes that the case $J=0$ 
corresponds to $a = b$ and $E = 2a^2 + \rmO(a^4)$ in the 
notation of \cite{GH1}.}
\[
  \lambda_1(E) \,=\, -E + \rmO(E^2)~, \quad
  \lambda_2(E) \,=\, 3E + \rmO(E^2)~, \quad
  \hbox{as } E \to 0~.
\]
\end{Lemma}

\begin{Remark}\label{th:GL}
By construction the wave profile $W_{J,E}(x) = \rme^{\rmi px}
Q_{J,E}(2kx)$ is a stationary solution of the time-dependent
Ginzburg-Landau equation
\begin{equation}\label{e:tdgl}
  W_t(x,t) \,=\, W_{xx}(x,t) + W(x,t) -|W(x,t)|^2 W(x,t)~.
\end{equation}
If we look for solutions of \reff{e:tdgl} of the form $W(x,t) =
\rme^{\rmi px} Q(2kx,t)$, we obtain the evolution equation
\begin{equation}\label{e:GLQ}
  Q_t(z,t) \,=\, 4k^2 Q_{zz}(z,t) + 4\rmi pk Q_z(z,t) +
  (1-p^2)Q(z,t) - |Q(z,t)|^2 Q(z,t)~,
\end{equation}
which is very similar to \reff{e:Qeq}. In particular, the
linearization of \reff{e:GLQ} at the equilibrium $Q_{J,E}$ is $Q_t
= -H_{J,E}Q$. Thus Proposition~\ref{th:Hneg} implies that the
quasi-periodic solutions of the Ginzburg-Landau equation
\reff{e:tdgl} are always {\em unstable}, even within the class of
solutions with the same period and the same Floquet multiplier as
the original wave profile. This complements previous results by
Bridges \& Rowlands \cite{BR1}, and by Doelman, Gardner \& Jones
\cite{DGJ}, which show that the quasi-periodic solutions of
\reff{e:tdgl} are unstable with respect to long-wave disturbances
(``sideband instability'').
\end{Remark}

We now continue with the proof of Proposition~\ref{th:Qorbit}.
The next observation is that, for any $(J,E) \in D$, the equilibrium
$Q_{J,E}$ of \reff{e:Qeq} is a member of a {\em two-parameter family}
of travelling and rotating waves of the form
\begin{equation}\label{e:travrot}
  Q(z,t) \,=\, \rme^{-\rmi \omega t} \,Q_{J,E}^{\omega,c}(z+ct)~,
  \quad z \in \R~, \quad t \in \R~,
\end{equation}
where $(\omega,c)$ lies in a neighborhood of the origin in $\R^2$
(the Lie algebra of $G$), and the profile $Q_{J,E}^{\omega,c} \in X$
is a smooth function of $(\omega,c)$ with $Q_{J,E}^{0,0} = Q_{J,E}$.
(Actually we even have a four-parameter family of such waves if
we take into account the action of the symmetry group $G$.)
Indeed, take $(J',E') \in D$ close to $(J,E)$, and define
$k', \ell'$ by the formulas \reff{e:kldef} with $(J,E)$ replaced by
$(J',E')$. Then
\[
  U(x,t) \,=\, \rme^{\rmi (p'x-t)} \,Q_{J'\!,E'}(2k'x)~, \quad
  \hbox{where } p' = k' + \ell'~,
\]
is a solution of \reff{e:nls}, but it is not of the form
\reff{e:Qgen} because $p' \neq p$ and $k' \neq k$ in general.
However we can transform $U(x,t)$ into a solution of \reff{e:nls}
of the form \reff{e:Qgen}, \reff{e:travrot} by applying successively
a dilation of factor $\lambda$ and a Galilean transformation of speed
$v$, where
\begin{equation}\label{e:lamvdef}
  \lambda \,=\, \frac{k}{k'}~, \quad v \,=\, 2(\lambda p'-p)~.
\end{equation}
After some elementary algebra, we obtain $Q_{J,E}^{\omega,c}(z) =
\lambda Q_{J'\!,E'}(z)$ with
\begin{equation}\label{e:omcdef}
  \omega \,=\, \lambda^2(1 - p'^2) - (1 - p^2)~, \quad
  c \,=\, 4 \lambda^2 k'p' - 4 kp~.
\end{equation}

The mapping $(J',E') \to (\omega,c)$ that we have just defined is a
diffeomorphism from a neighborhood of $(J,E)$ onto a neighborhood
of $(0,0)$. Indeed, using \reff{e:kldef} and \reff{e:lamvdef},
the formulas \reff{e:omcdef} can be written more explicitly as
\begin{equation}\label{e:omcexp}
  \omega \,=\, \frac{k^2}{\pi^2}\Bigl(T'^2 - (\pi+\Psi')^2
  \Bigr) - (1 - p^2)~, \quad c \,=\, \frac{4k^2}{\pi}(\pi+\Psi')
  - 4 kp~,
\end{equation}
where $T' = T(J',E')$ and $\Psi' = \Psi(J',E')$. Differentiating
these expressions with respect to $J',E'$ we obtain
\begin{equation}\label{e:MMdef}
  \MM_{J,E} \,\eqdef\, \pmatrix{
    \frac{\partial \omega}{\partial E'} &
    \frac{\partial c}{\partial E'} \vspace{1mm}\cr
    \frac{\partial \omega}{\partial J'} &
    \frac{\partial c}{\partial J'}}
  \bigg|_{(J',E') = (J,E)} ~=~ \pmatrix{\frac{2k}{\pi}
    (\frac{\partial T}{\partial E} -p\frac{\partial \Psi}{\partial E})&
    \frac{4k^2}{\pi}\frac{\partial \Psi}{\partial E} \vspace{1mm} \cr
    \frac{2k}{\pi}(\frac{\partial T}{\partial J} -p\frac{\partial \Psi}
    {\partial J}) & \frac{4k^2}{\pi}\frac{\partial \Psi}{\partial J}}~.
\end{equation}
By Proposition~\ref{th:KAM}, we have $\det(\MM_{J,E}) = (8k^3/\pi^2)
\Delta(J,E) > 0$, hence the matrix $\MM_{J,E}$ is invertible.
This proves the existence of the travelling and rotating waves
\reff{e:travrot} for sufficiently small $(\omega,c) \in \R^2$.

By construction, the profile $Q_{J,E}^{\omega,c} \in X$ is a critical
point of the functional
\begin{equation}\label{e:EJEomc}
  \EE_{J,E}^{\omega,c}(Q) \,=\, \EE_{J,E}(Q) - \omega N(Q) - c M(Q)~,
  \quad Q \in X~.
\end{equation}
For later use, we remark that, if $\lambda$ is as in \reff{e:lamvdef},
we have the identity
\begin{equation}\label{e:primeid}
  \EE_{J,E}^{\omega,c}(\lambda Q) \,=\, \lambda^4 \EE_{J'\!,E'}(Q)~,
  \quad \hbox{for all } Q\in X~.
\end{equation}
Following \cite{GSS2} we now define, for $(\omega,c)$ in a
neighborhood of $(0,0)$,
\begin{equation}\label{e:ddef}
  d_{J,E}(\omega,c) \,=\, \EE_{J,E}^{\omega,c}(Q_{J,E}^{\omega,c})
  \,\equiv\, \lambda^4 \EE_{J'\!,E'}(Q_{J'\!,E'})~.
\end{equation}
A crucial role in the orbital stability argument will be played by
the Hessian matrix of $d_{J,E}$, namely:
\[
   \HH_{J,E} \,=\, \pmatrix{
    \frac{\partial^2 d_{J,E}}{\partial \omega^2} &
    \frac{\partial^2 d_{J,E}}{\partial \omega \,\partial c} \vspace{2mm}\cr
    \frac{\partial^2 d_{J,E}}{\partial c \,\partial \omega} &
    \frac{\partial^2 d_{J,E}}{\partial c^2}}
  \Bigg|_{(\omega,c) = (0,0)}~.
\]

\begin{Proposition}\label{th:Hessian}
For all $(J,E) \in D$ we have $\det(\HH_{J,E}) < 0$.
\end{Proposition}

\begin{Remark}
Since $\HH_{J,E}$ is a symmetric matrix, Proposition~\ref{th:Hessian}
implies of course that $\HH_{J,E}$ has one positive and one negative
eigenvalue. This can be verified explicitly in the case of small
amplitude periodic waves, for which the following asymptotic result
is established in \cite{GH1}:
\[
  \HH_{J,E} \,=\, \frac{\pi}3 \pmatrix{-2 & -1 \cr -1 & 1}
  (\1 + \rmO(E))~, \quad \hbox{as} \quad (J,E) \to (0,0)
  \hbox{ in }D~.
\]
\end{Remark}

Since the proof of Proposition~\ref{th:Hessian} is long and technical,
we postpone it to the end of this section, and we now show how
Propositions~\ref{th:Hneg} and \ref{th:Hessian} together imply
Proposition~\ref{th:Qorbit}. The arguments here are rather classical
and can be found in \cite{GSS2}, so we shall just indicate how the
general theory of \cite{GSS2} applies to present case.

For any $(J,E) \in D$, Proposition~\ref{th:Hneg} shows that the wave
profile $Q_{J,E}$ is a degenerate saddle point of the energy
$\EE_{J,E}$, with one unstable and two neutral directions
(the latter are due to the fact that \reff{e:Qeq} is $G$-invariant).
To get rid of the unstable direction we remark that the
evolution of \reff{e:Qeq} does not take place in the whole function
space $X$, but on codimension two surfaces where the charge $N$
and the momentum $M$ are constant. Let
\begin{equation}\label{e:SJEdef}
  \Sigma_{J,E} \,=\, \Big\{Q \in X\,\Big|\, N(Q) = N(Q_{J,E})\,,~
  M(Q) = M(Q_{J,E})\Big\}~.
\end{equation}
It is easy to verify that $\Sigma_{J,E}$ is indeed a smooth
submanifold of $X$ of codimension two in a neighborhood of $Q_{J,E}$,
or more generally in a neighborhood of the orbit of $Q_{J,E}$ under
$G$. The crucial point is that the functional $\EE_{J,E}$ is coercive
on $\Sigma_{J,E}$ with respect to the semi-distance $\rho$ defined
by \reff{e:rhodef}. More precisely, there exist positive constants
$C_1, \delta$ (depending on $J,E$) such that
\begin{equation}\label{e:coercive}
  \EE_{J,E}(Q) - \EE_{J,E}(Q_{J,E}) \,\ge\, C_1 \rho(Q,Q_{J,E})^2~,
\end{equation}
for all $Q \in \Sigma_{J,E}$ such that $\rho(Q,Q_{J,E}) \le \delta$.

To prove \reff{e:coercive} we first note that the tangent space to
$\Sigma_{J,E}$ at $Q_{J,E}$ is
\[
  \TT_{J,E} \,=\, \Bigl\{Q \in X \,\Big|\, \langle Q_{J,E},Q\rangle
  = \langle \rmi Q_{J,E}',Q\rangle = 0\Bigr\}~,
\]
(see \reff{e:firstder}). We next introduce the ``normal space''
$\NN_{J,E} = \linspan\{\partial_\omega Q_{J,E},\partial_c
Q_{J,E}\}$, where
\[
  \partial_\omega Q_{J,E} \,=\, \frac{\partial}{\partial \omega}
  Q_{J,E}^{\omega,c}\Big|_{(\omega,c)=(0,0)}~, \quad
  \partial_c Q_{J,E} \,=\, \frac{\partial}{\partial c}
  Q_{J,E}^{\omega,c}\Big|_{(\omega,c)=(0,0)}~.
\]
As $Q_{J,E}^{\omega,c}$ is a critical point of $\EE_{J,E}^{\omega,c}$,
we have $\EE_{J,E}'(Q_{J,E}^{\omega,c}) = \omega
N'(Q_{J,E}^{\omega,c}) + c M'(Q_{J,E}^{\omega,c})$. Differentiating
this relation with respect to $\omega$ and $c$ at $(\omega,c) = (0,0)$,
we obtain
\[
  H_{J,E}(\partial_\omega Q_{J,E}) \,=\, N'(Q_{J,E}) \,=\, Q_{J,E}~, \quad
  H_{J,E}(\partial_c Q_{J,E}) \,=\, M'(Q_{J,E}) \,=\, \rmi
  Q_{J,E}'~,
\]
hence $\langle H_{J,E}Q_1,Q_2\rangle = 0$ for all $Q_1 \in \NN_{J,E}$,
$Q_2 \in \TT_{J,E}$. In a similar way, differentiating \reff{e:ddef}
we find
\begin{equation}\label{e:dNM}
  \frac{\partial}{\partial \omega}d_{J,E}(\omega,c) \,=\,
    -N(Q_{J,E}^{\omega,c})~, \quad
  \frac{\partial}{\partial c}d_{J,E}(\omega,c) \,=\,
    -M(Q_{J,E}^{\omega,c})~,
\end{equation}
hence
\begin{eqnarray}\label{e:Hexp1}
  \HH_{J,E} &=&  -\pmatrix{
    \frac{\partial}{\partial \omega} N(Q_{J,E}^{\omega,c}) &
    \frac{\partial}{\partial \omega} M(Q_{J,E}^{\omega,c})\vspace{1mm}\cr
    \frac{\partial}{\partial c} N(Q_{J,E}^{\omega,c}) &
    \frac{\partial}{\partial c} M(Q_{J,E}^{\omega,c})}
  \Bigg|_{(\omega,c) = (0,0)} \\[2mm] \label{e:Hexp2}
  &=& -\pmatrix{
  \langle H_{J,E}(\partial_\omega Q_{J,E}),\partial_\omega Q_{J,E}\rangle &
  \langle H_{J,E}(\partial_c Q_{J,E}),\partial_\omega Q_{J,E}\rangle \cr
  \langle H_{J,E}(\partial_\omega Q_{J,E}),\partial_c Q_{J,E}\rangle &
  \langle H_{J,E}(\partial_c Q_{J,E}),\partial_c Q_{J,E}\rangle}~.
\end{eqnarray}
Thus $-\HH_{J,E}$ is the matrix of the bilinear form $b_{J,E}(Q_1,Q_2)
= \langle H_{J,E} Q_1,Q_2\rangle$ restricted to the subspace
$\NN_{J,E}$. As $\HH_{J,E}$ is non-degenerate by
Proposition~\ref{th:Hessian}, we have $\NN_{J,E} \cap \TT_{J,E} =
\{0\}$, hence also $X = \NN_{J,E} \oplus \TT_{J,E}$.  Moreover, since
$H_{J,E}$ has only one negative eigenvalue by
Proposition~\ref{th:Hneg}, and since the restriction of $b_{J,E}$ to
$\NN_{J,E}$ {\em has} a negative eigenvalue by
Proposition~\ref{th:Hessian}, a standard signature argument shows that
the restriction of $b_{J,E}$ to the complementary space $\TT_{J,E}$ is
nonnegative. More precisely, if we decompose $\TT_{J,E} = K_{J,E}
\oplus \TT_{J,E}'$ where $K_{J,E} = \linspan\{Q_{J,E}',\rmi Q_{J,E}\}$
is the kernel of $H_{J,E}$ in $X$ and $\TT_{J,E}' = \{Q \in \TT_{J,E}
\,|\, \langle Q_{J,E}',Q\rangle = \langle \rmi Q_{J,E},Q\rangle =
0\}$, then there exists $C_2 > 0$ such that
\begin{equation}\label{e:coercive2}
  \langle H_{J,E} Q,Q\rangle \,\ge\, C_2 \|Q\|_X^2~, \quad
  \hbox{for all } Q \in \TT_{J,E}'~.
\end{equation}
Estimate \reff{e:coercive} follows from its infinitesimal version
\reff{e:coercive2} using a Taylor expansion of the energy
$\EE_{J,E}$ at $Q_{J,E}$ and the fact that $\EE_{J,E}$ is
$G$-invariant (see \cite{GSS2} for details).

\begin{Remark}\label{th:C1C2}
The argument above does not give any lower bound on the constants
$C_1, C_2$ in \reff{e:coercive}, \reff{e:coercive2}. In particular,
we cannot take for $C_2$ the lowest positive eigenvalue of the
operator $H_{J,E}$, unless the subspace $\NN_{J,E}$ contains precisely
the eigenfunction associated to the (unique) negative eigenvalue
of $H_{J,E}$. However, one can check that the constants $C_1, \delta$
in \reff{e:coercive} are bounded away from zero uniformly for $(J,E)$
in any compact subdomain of $D$.
\end{Remark}

With inequality \reff{e:coercive} at hand, it is now easy to conclude
the proof of Proposition~\ref{th:Qorbit}. Assume that $Q_0 \in X$
satisfies $\rho(Q_0,Q_{J,E}) \le \epsilon$ for some small
$\epsilon > 0$. Replacing $Q_0$ by $\RR_{(\phi,\xi)}Q_0$ if needed,
we can assume that $\|Q_0-Q_{J,E}\|_X \le \epsilon$. As $Q_{J,E}$ is
a critical point of $\EE_{J,E}$, we have $|\EE_{J,E}(Q_0) -
\EE_{J,E}(Q_{J,E})| \le C_3 \epsilon^2$, for some $C_3 > 0$
(depending on $J,E$). We distinguish two cases:\\[1mm]
{\bf 1.} If $Q_0 \in \Sigma_{J,E}$, then the solution of \reff{e:Qeq}
with initial data $Q_0$ satisfies $Q(\cdot,t) \in \Sigma_{J,E}$ for
all $t \in \R$, and \reff{e:coercive} implies that
\[
  C_1 \rho(Q(\cdot,t),Q_{J,E})^2 \,\le\,
  \EE_{J,E}(Q(\cdot,t)) - \EE_{J,E}(Q_{J,E})
  \,=\, \EE_{J,E}(Q_0) - \EE_{J,E}(Q_{J,E}) \,\le\, C_3 \epsilon^2~,
\]
for all $t \in\R$, provided $\epsilon$ is small enough so that
$C_3\epsilon^2 < C_1 \delta^2$. This is the desired result.
\\[1mm]
{\bf 2.} The case where $Q_0 \notin \Sigma_{J,E}$ can be reduced
to the previous one by the following argument. In view of
\reff{e:Hexp1} and Proposition~\ref{th:Hessian}, the map $(\omega,c)
\mapsto (N(Q_{J,E}^{\omega,c}),M(Q_{J,E}^{\omega,c}))$ is a diffeomorphism
from a neighborhood of $(0,0)$ onto a neighborhood of $(N(Q_{J,E}),
M(Q_{J,E}))$. Thus there exists a unique $(\omega,c) \in \R^2$
with $|\omega|+|c| = \rmO(\epsilon)$ such that $N(Q_0) =
N(Q_{J,E}^{\omega,c})$, $M(Q_0) = M(Q_{J,E}^{\omega,c})$. By
construction $Q_{J,E}^{\omega,c} = \lambda Q_{J'\!,E'}$, where
$(J',E')$ is $\epsilon$-close to $(J,E)$ and $\lambda$ is as in
\reff{e:lamvdef}. Thus, if we denote $\tilde Q(z,t) = \lambda^{-1}
Q(z,t)$, we see that $\tilde Q(\cdot,t) \in \Sigma_{J'\!,E'}$
for all $t \in \R$. Using \reff{e:primeid} together with
\reff{e:coercive}, we thus obtain
\begin{eqnarray*}
  \EE_{J,E}^{\omega,c}(Q(\cdot,t)) - \EE_{J,E}^{\omega,c}(Q_{J,E}^{\omega,c})
  &=& \lambda^4 \Bigl(\EE_{J'\!,E'}(\tilde Q(\cdot,t)) -
  \EE_{J'\!,E'}(Q_{J'\!,E'})\Bigr) \\
  &\ge& C_1' \lambda^4  \rho(\tilde Q(\cdot,t),Q_{J'\!,E'})^2
  ~=~ C_1' \lambda^2 \rho(Q(\cdot,t),Q_{J,E}^{\omega,c})^2~,
\end{eqnarray*}
as long as $\rho(Q(\cdot,t),Q_{J,E}^{\omega,c})$ stays sufficiently
small. Since
\[
  \EE_{J,E}^{\omega,c}(Q(\cdot,t)) - \EE_{J,E}^{\omega,c}(Q_{J,E}^{\omega,c})
  \,=\, \EE_{J,E}^{\omega,c}(Q_0) - \EE_{J,E}^{\omega,c}(Q_{J,E}^{\omega,c})
  \,\le\, C_3'\|Q_0 - Q_{J,E}^{\omega,c}\|_X^2~,
\]
and $\|Q_0 - Q_{J,E}^{\omega,c}\|_X \le \|Q_0 - Q_{J,E}\|_X +
\|Q_{J,E} - Q_{J,E}^{\omega,c}\|_X \le C_4 \epsilon$, we conclude
that
\[
  \rho(Q(\cdot,t),Q_{J,E}) \,\le\,
  \rho(Q(\cdot,t),Q_{J,E}^{\omega,c}) +
  \|Q_{J,E}^{\omega,c} - Q_{J,E}\|_X \,\le\, C_5 \epsilon
  \quad \hbox{for all } t \in\R~,
\]
if $\epsilon > 0$ is sufficiently small. This concludes the proof
of Proposition~\ref{th:Qorbit}. \QED

\medskip
\noindent{\bf Proof of Proposition~\ref{th:Hessian}.}
We start from the expression \reff{e:Hexp1} of $\HH_{J,E}$.
We recall that $(\omega,c)$ can be parametrized (in a neighborhood
of the origin) by $(J',E')$ according to the formulas \reff{e:omcdef}
or \reff{e:omcexp}, and that $Q_{J,E}^{\omega,c} = \lambda
Q_{J'\!,E'}$ where $\lambda$ is given by \reff{e:lamvdef}.
Thus $\HH_{J,E} = -(\MM_{J,E})^{-1}\KK_{J,E}$, where $\MM_{J,E}$
is defined in \reff{e:MMdef} and
\[
  \KK_{J,E} \,=\, \pmatrix{
   \frac{\partial' N}{\partial' E} &
   \frac{\partial' M}{\partial' E} \vspace{1mm}\cr
   \frac{\partial' N}{\partial' J} &
   \frac{\partial' M}{\partial' J}} ~\eqdef~ \pmatrix{
   \frac{\partial}{\partial E'} N(\lambda Q_{J'\!,E'}) &
   \frac{\partial}{\partial E'} M(\lambda Q_{J'\!,E'}) \vspace{1mm} \cr
   \frac{\partial}{\partial J'} N(\lambda Q_{J'\!,E'}) &
   \frac{\partial}{\partial J'} M(\lambda Q_{J'\!,E'})}
   \bigg|_{(J',E') = (J,E)}~.
\]
Since we already know that $\det(\MM_{J,E}) > 0$, it remains to verify
that $\det(\KK_{J,E}) < 0$.

We first give more explicit formulas for the coefficients of
$\KK_{J,E}$. As $W_{J,E}(x) = \rme^{\rmi px}Q_{J,E}(2kx)$, where
$p = k+\ell$ and $k,\ell$ are given by \reff{e:kldef}, it follows
from \reff{e:NMEdef} that
\begin{eqnarray}\nonumber
  N(Q_{J,E}) &=& \frac12 \int_0^{2\pi} |Q_{J,E}(z)|^2\dd z
  \,=\, k \int_0^T |W_{J,E}(x)|^2 \dd x~, \\\nonumber
  M(Q_{J,E}) &=& \frac{\rmi}2 \int_0^{2\pi} \overline{Q}_{J,E}(z)
  Q_{J,E}'(z)\dd z \\\label{e:Mexp}
  &=& \frac{\rmi}2 \int_0^T \Bigl(\overline{W}_{J,E}(x)W_{J,E}'(x)
  -\rmi p |W_{J,E}(x)|^2\Bigr)\dd x \,=\, \frac{p}{2k}N(Q_{J,E})
  - \frac12 JT~.
\end{eqnarray}
Moreover, proceeding as in Lemma~\ref{th:TPhi} we find
\begin{equation}\label{e:Nexp}
  N(Q_{J,E}) \,=\, 2k \int_{r_1(J,E)}^{r_2(J,E)} \frac{r^2}
  {\sqrt{2(E - V_J(r))}}\dd r \,=\, 2\sqrt{2}\,k \int_0^{\pi/2}
  \frac{s(\phi)}{\sqrt{y_3 - s(\phi)}}\,\dd \phi~,
\end{equation}
where $s(\phi) = y_1\cos^2(\phi) + y_2\sin^2(\phi)$ and $y_1, y_2,
y_3$ are as in Fig.~3. Of course, since $N,M$ are quadratic
functionals, we have
\begin{equation}\label{e:NMexp}
  N(\lambda Q_{J'\!,E'}) \,=\, \lambda^2 N(Q_{J'\!,E'})
  \,=\, \frac{k^2}{k'^2}\,N(Q_{J'\!,E'})~, \quad
  M(\lambda Q_{J'\!,E'}) \,=\, \frac{k^2}{k'^2}\,M(Q_{J'\!,E'})~.
\end{equation}
Differentiating \reff{e:NMexp} with respect to $J',E'$ and
using \reff{e:Mexp}, \reff{e:Nexp}, we obtain after straightforward
calculations:
\begin{eqnarray}\label{e:Nder}
  \frac{\partial' N}{\partial' E} &=&
  C_1 \frac{\partial y_1}{\partial E} +
  C_2 \frac{\partial y_2}{\partial E} +
  C_3 \frac{\partial T}{\partial E}~, \quad
  \frac{\partial' N}{\partial' J} ~=~
  C_1 \frac{\partial y_1}{\partial J} +
  C_2 \frac{\partial y_2}{\partial J} +
  C_3 \frac{\partial T}{\partial J}~, \\[2mm]\label{e:Mder}
  \frac{\partial' M}{\partial' E} &=&
  \frac{p}{2k}\frac{\partial' N}{\partial' E} +
  \frac{N}{2\pi}\frac{\partial \Psi}{\partial E} -
  \frac{3J}2 \frac{\partial T}{\partial E}~, \quad
  \frac{\partial' M}{\partial' J} ~=~
  \frac{p}{2k}\frac{\partial' N}{\partial' J} +
  \frac{N}{2\pi}\frac{\partial \Psi}{\partial J} -
  \frac{3J}2 \frac{\partial T}{\partial J} - \frac{T}2~,
\end{eqnarray}
where $N = N(Q_{J,E})$ and
\begin{eqnarray*}
  C_1 &=& \sqrt{2}\,k \int_0^{\pi/2} \frac{2y_3 \cos^2(\phi) +
  s(\phi)\sin^2(\phi)}{(y_3-s(\phi))^{3/2}}\dd\phi~, \\
  C_2 &=& \sqrt{2}\,k \int_0^{\pi/2} \frac{2y_3 \sin^2(\phi) +
  s(\phi)\cos^2(\phi)}{(y_3-s(\phi))^{3/2}}\dd\phi~, \\
  C_3 &=& 2\sqrt{2} \,\frac{k^2}{\pi} \int_0^{\pi/2} \frac{s(\phi)}
  {\sqrt{y_3-s(\phi)}}\dd\phi \,=\, \frac{k}{\pi}\,N(Q_{J,E})~.
\end{eqnarray*}

We now compute $\det(\KK_{J,E})$. Replacing the second column of the
matrix $\KK_{J,E}$ by its expression \reff{e:Mder}, we see that
$\det(\KK_{J,E}) = \Delta_1 + \Delta_2 + \Delta_3$, where
\[
  \Delta_1 \,=\, \Bigl(\frac{NB_3}{2\pi} -\frac{T}2\Bigr)
  \frac{\partial' N}{\partial' E}~,
\]
and
\[
  \Delta_2 \,=\,\frac{N}{2\pi}\det\pmatrix{
  \frac{\partial' N}{\partial' E} &
  \frac{\partial \Psi}{\partial E} \vspace{1mm} \cr
  \frac{\partial' N}{\partial' J} &
  \frac{\partial \Psi}{\partial J} - B_3}~, \quad
  \Delta_3 \,=\,-\frac{3J}{2}\det\pmatrix{
  \frac{\partial' N}{\partial' E} &
  \frac{\partial T}{\partial E} \vspace{1mm} \cr
  \frac{\partial' N}{\partial' J} &
  \frac{\partial T}{\partial J}}~.
\]
Here $B_3$ is as in \reff{e:Psider}.

We claim that $\Delta_1 < 0$. Indeed, arguing as in \reff{e:A1A2}
and using Lemma~\ref{th:fg}, we find
\[
  C_2 - C_1 \,=\, \sqrt{2}\,k \int_0^{\pi/2}\frac{2y_3-s(\phi)}
  {(y_3-s(\phi))^{3/2}}\,(\sin^2(\phi)-\cos^2(\phi))\dd\phi \,>\, 0~,
\]
because $\phi \mapsto (2y_3-s(\phi))(y_3-s(\phi))^{-3/2}$ is strictly
increasing over $[0,\pi/2]$. Thus $C_2 > C_1 > 0$, and proceeding
as in \reff{e:TEpos} we deduce that $C_1(\partial y_1/\partial E)
+ C_2(\partial y_2/\partial E) > 0$. Since $C_3 > 0$ and $\partial T/
\partial E > 0$ by Proposition~\ref{th:monotone}, we conclude that
$\partial' N/\partial' E > 0$. On the other hand, since $s(\phi) \le
y_2$ for all $\phi\in [0,2\pi]$, it follows from \reff{e:Texp},
\reff{e:kldef}, \reff{e:Nexp} that $N \le \pi y_2$. Moreover, as
\[
  \frac{y_2}{\sqrt{y_3}(\sqrt{y_3} + \sqrt{y_3-s(\phi)})} \,\le\,
  \frac{y_2}{\sqrt{y_3}(\sqrt{y_3} + \sqrt{y_3-y_2})} \,=\,
  \frac{\sqrt{y_3}-\sqrt{y_3-y_2}}{\sqrt{y_3}} \,<\, 1~,
\]
it follows from \reff{e:Texp}, \reff{e:Bidef} that $B_3 y_2 < T$.
Thus $NB_3 < \pi T$, hence $\Delta_1 < 0$.

It remains to verify that $\Delta_2 + \Delta_3 < 0$. By \reff{e:Tder},
\reff{e:Nder} we have
\[
  \frac{\partial' N}{\partial' E} \,=\,
  (C_1+C_3 A_1) \frac{\partial y_1}{\partial E} +
  (C_2+C_3 A_2) \frac{\partial y_2}{\partial E}~,
\]
and similarly for the $J$-derivative. Thus using \reff{e:Psider}
and proceeding as in \reff{e:Deltaexp} we find
\[
  \Delta_2 \,=\, \frac{N}{2\pi}\,(y_2-y_1)\frac{\partial y_1}{\partial J}
  \frac{\partial y_2}{\partial J}\Bigl((C_1+C_3 A_1)B_2 -
  (C_2+C_3A_2) B_1\Bigr)~.
\]
In a similar way
\[
  \Delta_3 \,=\, -\frac32\,(y_2-y_1)\frac{\partial y_1}{\partial J}
  \frac{\partial y_2}{\partial J}(C_1 A_2 - C_2 A_1)~.
\]
Thus $\Delta_2 + \Delta_3 = (y_2-y_1)(\partial y_1/\partial J)
(\partial y_2/\partial J)\Delta_4$, where
\[
  \Delta_4 \,=\, \frac32 (A_1 C_2 - A_2 C_1) + \frac{NC_3}{2\pi}
  (A_1 B_2 - A_2 B_1) + \frac{N}{2\pi} (C_1 B_2 - C_2 B_1)~.
\]
As $y_2 > y_1$ and $(\partial y_1/\partial J)(\partial y_2/\partial J)
< 0$, we have to verify that $\Delta_4 > 0$.

Like in the proof of Proposition~\ref{th:KAM}, it is convenient to
express the various constants $A_i, B_i, C_i$ in terms of the new
function $\sigma(\phi) = (1-y_3^{-1}s(\phi))^{1/2}$. Then $A_i =
\sqrt{2}\,y_3^{-3/2}\bA_i$, $B_i = \sqrt{2}\,y_3^{-5/2}\bB_i$,
$C_i = \sqrt{2}\,k y_3^{-1/2}\bC_i$, where $\bA_i$, $\bB_i$ are
defined in \reff{e:bAiBidef} and
\begin{eqnarray*}
  \bC_1 &=& \int_0^{\pi/2} \frac{2\cos^2(\phi) + (1{-}\sigma(\phi)^2)
  \sin^2(\phi)}{\sigma(\phi)^3}\,\dd\phi~, \\
  \bC_2 &=& \int_0^{\pi/2} \frac{2\sin^2(\phi) + (1{-}\sigma(\phi)^2)
  \cos^2(\phi)}{\sigma(\phi)^3}\,\dd\phi~.
\end{eqnarray*}
Thus
\[
  \Delta_4 \,=\, \frac{2k}{y_3^2}\left(\frac32
  (\bA_1 \bC_2 - \bA_2 \bC_1) + \frac{NC_3}{2\pi k y_3^2}
  (\bA_1 \bB_2 - \bA_2 \bB_1) + \frac{N}{2\pi y_3}
  (\bC_1 \bB_2 - \bC_2 \bB_1)\right)~.
\]
Unfortunately, the three terms in the right-hand side are not all
positive. Indeed, we recall that $\bB_i = \bA_i + \tilde \bB_i$
for $i = 1,2$, where $\tilde \bB_1, \tilde \bB_2$ are defined in
\reff{e:tildeBidef}. Similarly, $\bC_i = \bA_i - \tilde \bC_i$, where
\begin{equation}\label{e:tildeCidef}
  \tilde \bC_1 \,=\,\int_0^{\pi/2} \frac{\sin^2(\phi)}{\sigma(\phi)}
  \,\dd\phi~, \quad
  \tilde \bC_2 \,=\,\int_0^{\pi/2} \frac{\cos^2(\phi)}{\sigma(\phi)}
  \,\dd\phi~.
\end{equation}
We also recall that $\bA_2 > \bA_1$, $\bB_2 > \bB_1$, $\bC_2 > \bC_1$,
but $\tilde \bB_1 > \tilde \bB_2$, $\tilde \bC_1 > \tilde \bC_2$.
All these bounds follow from Lemma~\ref{th:fg}, as in
\reff{e:A1A2}. Thus
\begin{eqnarray}\nonumber
  \bA_1 \bC_2 - \bA_2 \bC_1 &=& \tilde \bC_1 \bA_2 - \tilde \bC_2
  \bA_1 \,>\, 0~, \\\label{e:ABCineq}
  \bA_1 \bB_2 - \bA_2 \bB_1 &=& \bA_1 \tilde \bB_2 - \bA_2 \tilde
  \bB_1 \,<\, 0~, \\\nonumber
  \bC_1 \bB_2 - \bC_2 \bB_1 &=& (\bA_1 \bB_2 - \bA_2 \bB_1) +
  (\tilde \bC_2 \bB_1 - \tilde \bC_1 \bB_2) \,<\, 0~.
\end{eqnarray}
As $N \le \pi y_2 < \pi y_3$ and $C_3 = (kN/\pi) < ky_3$, it follows
from \reff{e:ABCineq} that $\Delta_4 > (2k/y_3^2)\Delta_5$ where
\[
  \Delta_5 \,=\, \frac32 (\bA_1 \bC_2 - \bA_2 \bC_1) +
  \frac12 (\bA_1 \bB_2 - \bA_2 \bB_1) +
  \frac12 (\bC_1 \bB_2 - \bC_2 \bB_1)~.
\]

It remains to verify that $\Delta_5 > 0$. Since $\bB_i = \bA_i + \tilde \bB_i$
and $\bC_i = \bA_i - \tilde \bC_i$ for $i = 1,2$, we have
\begin{eqnarray*}
  \Delta_5 &=&  (\tilde \bC_1 \bA_2  - \tilde \bC_2 \bA_1) +
  (\bA_1 \tilde \bB_2 - \bA_2 \tilde \bB_1) +
  \frac12 (\tilde \bB_1 \tilde \bC_2 - \tilde \bB_2 \tilde \bC_1) \\
  &=& \Bigl(\bA_2 - \bA_1 + \frac{\tilde \bC_1-\tilde \bC_2}2\Bigr)
  (\tilde \bC_1-\tilde \bB_1) + \Bigl(\bA_1 - \frac{\tilde \bC_1}2
  \Bigr)\Bigl((\tilde \bC_1 - \tilde \bC_2) -
  (\tilde \bB_1 - \tilde \bB_2)\Bigr)~.
\end{eqnarray*}
We claim that all terms in the right-hand side are now positive.
Indeed, we already know that $\bA_2 - \bA_1 > 0$ and $\tilde \bC_1 -
\tilde \bC_2 > 0$. Using \reff{e:tildeBidef} and \reff{e:tildeCidef}
we obtain
\[
  \tilde \bC_1-\tilde \bB_1 \,=\, \int_0^{\pi/2}\frac{(2+\sigma(\phi))
  \sin^2(\phi)-1}{(1+\sigma(\phi))^2}\dd \phi \,>\,
  \int_0^{\pi/2}\frac{\sin^2(\phi) - \cos^2(\phi)}{(1+\sigma(\phi))^2}
  \dd \phi \,>\, 0~,
\]
because $\phi \mapsto (1+\sigma(\phi))^{-2}$ is strictly increasing
on $[0,\pi/2]$. On the other hand, since $\sigma(\phi) \le 1$ we have
\[
  \bA_1 \,>\, \int_0^{\pi/2}\frac{1}{\sigma(\phi)^3}\dd \phi \,>\,
  \int_0^{\pi/2}\frac{\sin^2(\phi)}{\sigma(\phi)}\dd \phi \,=\,
  \tilde \bC_1~,
\]
hence $\bA_1 - \tilde \bC_1/2 > 0$ a fortiori. Finally,
\[
  (\tilde \bC_1-\tilde \bC_2) -  (\tilde \bB_1-\tilde \bB_2)
  \,=\, \int_0^{\pi/2}\frac{2+\sigma(\phi)}{(1+\sigma(\phi))^2}
  \,(\sin^2(\phi) - \cos^2(\phi))\dd \phi \,>\, 0~,
\]
because $\phi \mapsto (2+\sigma(\phi))/(1+\sigma(\phi))^2$
is strictly increasing on $[0,\pi/2]$. Thus $\Delta_5 > 0$, and
the proof of Proposition~\ref{th:Hessian} (hence of
Proposition~\ref{th:Qorbit}) is now complete. \QED

\section{The focusing NLS case}
\label{s:fe}

In this final section we show how the preceding results can be
extended to the focusing NLS equation \reff{e:nls} with
$\gamma = 1$. In this case, Eq.\reff{e:snls} has quasi-periodic
solutions for all values of $\omega$, but we restrict ourselves
to the generic cases $\omega = 1$ and $\omega = -1$, which we
consider separately.

\subsection{Counter-rotating waves ($\omega = 1$)}
\label{s:om1}

Proceeding as in Section~\ref{s:ex}, we first study the solutions
of the stationary equation
\begin{equation}\label{e:sgl1}
  W_{xx}(x) + W(x) + |W(x)|^2 W(x) \,=\, 0~, \quad x \in \R~.
\end{equation}
This is again an integrable Hamiltonian system with conserved
quantities
\begin{equation}\label{e:EJ1}
  J \,=\, \Im(\overline{W}W_x)~, \qquad
  E \,=\, \frac12|W_x|^2 + \frac12|W|^2 + \frac14|W|^4~.
\end{equation}
In particular, the effective potential $V_J(r) = J^2/(2r^2) +
r^2/2 + r^4/4$ is now strictly convex for any $J \in \R$. If we
set $J = q(q^2-1)$ where $q \in \R$, $|q| \ge 1$, then the unique
minimum of $V_J$ is attained at $r = r_q = \sqrt{q^2-1}$ and has the value
\[
  E_-(J) \,=\, V_J(\sqrt{q^2-1}) \,=\, \frac14 (q^2-1)(3q^2+1)~.
\]
It follows that \reff{e:sgl1} has quasi-periodic solutions
if and only if $(J,E) \in D$, where
\[
  D \,=\, \Bigl\{(J,E) \in \R^2 \,\Big|\, E > E_-(J)\Bigr\}~.
\]
The period $T$ and the phase increment $\Phi$ of these solutions
are given by
\begin{equation}\label{e:TPhiexp1}
  T(J,E) \,=\, 2\sqrt{2} \int_0^{\pi/2}\frac{\dd \phi}{\sqrt{s(\phi)
  -y_3}}~, \quad
  \Phi(J,E) \,=\, 2\sqrt{2} \int_0^{\pi/2} \frac{J}{s(\phi)}\,
  \frac{\dd \phi}{\sqrt{s(\phi)-y_3}}~,
\end{equation}
where $s(\phi) = y_1 \cos^2(\phi) + y_2 \sin^2(\phi)$ and
$y_3 < 0 \le y_1 < y_2$ are the roots of the cubic polynomial
$P(y) = -y^3 -2y^2 + 4Ey - 2J^2$. Similarly, the renormalized
phase \reff{e:Psidef} satisfies
\begin{equation}\label{e:Psiexp1}
  \Psi(J,E) \,=\, -2\sqrt{2} \,J \int_0^{\pi/2} \frac{\dd \phi}
  {\sqrt{-y_3}\,\sqrt{s(\phi)-y_3}\,(\sqrt{-y_3} + \sqrt{s(\phi)-y_3})}~.
\end{equation}

In contrast with the defocusing case, the period $T(J,E)$ is now
a {\em decreasing} function of the energy. The analogue of
Proposition~\ref{th:monotone} is:

\begin{Proposition}\label{th:monot1}$ $\\[1mm]
i) $~\displaystyle \frac{\partial T}{\partial E}(J,E) < 0~$ for all
$(J,E) \in D$. \\[1mm]
ii) $~\displaystyle \frac{\partial \Psi}{\partial E}(J,E) =
-\frac{\partial T}{\partial J}(J,E) > 0~$ for all $(J,E) \in D$ with
$J > 0$.
\end{Proposition}

\figurewithtex 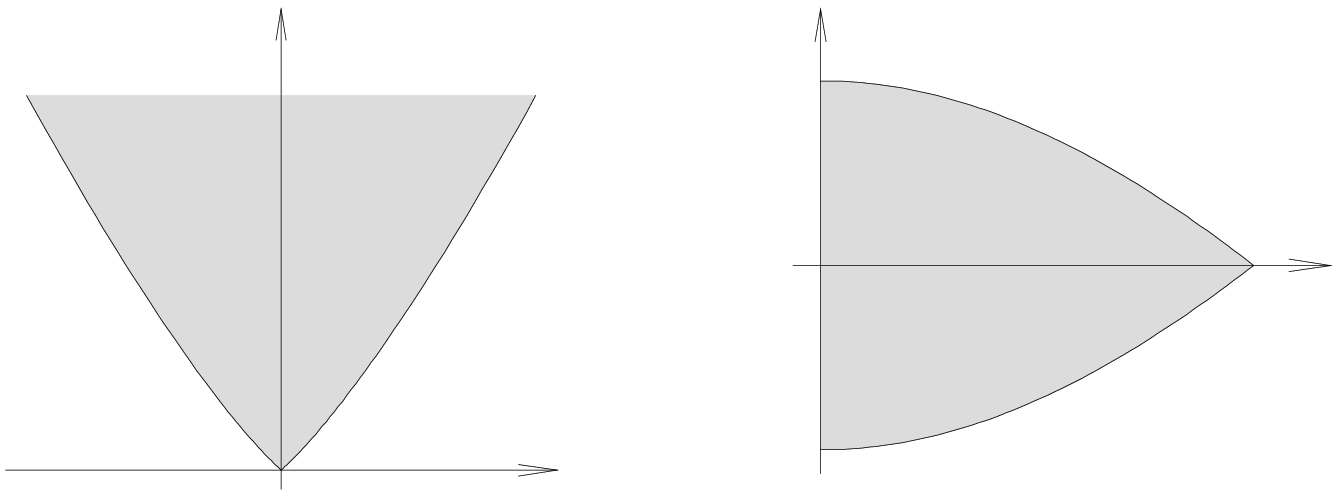 Fig5.tex 5.000 14.000
{\bf Fig.~5:} Existence domain for the counter-rotating waves of
the focusing NLS equation, in terms of the parameters $J,E$ (left)
and $T,\Psi$ (right).\cr

\noindent{\bf Proof.} Differentiating $T$ with respect to $E$
and $J$ we obtain
\begin{equation}\label{e:Tder1}
  \frac{\partial T}{\partial E} \,=\, -A_1 \frac{\partial y_1}{\partial E}
   -A_2 \frac{\partial y_2}{\partial E}~, \quad
  \frac{\partial T}{\partial J} \,=\, -A_1 \frac{\partial y_1}{\partial J}
  - A_2 \frac{\partial y_2}{\partial J}~, \quad
\end{equation}
where
\begin{equation}\label{e:Aiexp1}
  A_1 \,=\, \sqrt{2}\int_0^{\pi/2} \frac{1+\cos^2(\phi)}{(s(\phi)-y_3)^{3/2}}
  \dd\phi~, \quad
  A_2 \,=\, \sqrt{2}\int_0^{\pi/2} \frac{1+\sin^2(\phi)}{(s(\phi)-y_3)^{3/2}}
  \dd\phi~. \quad
\end{equation}
We remark that $\phi \mapsto (s(\phi)-y_3)^{-3/2}$ is strictly
decreasing over $[0,\pi/2]$. As
\begin{equation}\label{e:A1A21}
  A_1 - A_2 \,=\, \sqrt{2}\int_0^{\pi/2}\frac{\cos^2(\phi)-\sin^2(\phi)}
  {(s(\phi)-y_3)^{3/2}}\dd\phi~, \quad
\end{equation}
it follows from Lemma~\ref{th:fg} that $A_1 > A_2 > 0$. On the other
hand, since
\begin{equation}\label{e:yiEJder}
  \frac{\partial y_i}{\partial E} \,=\, -\frac{4y_i}{-3y_i^2 - 4y_i
  + 4E}~, \quad
  \frac{\partial y_i}{\partial J} \,=\, \frac{4J}{-3y_i^2 - 4y_i
  + 4E}~, \quad i = 1,2,3~,
\end{equation}
it is not difficult to verify that \reff{e:yEsign} still holds,
whereas \reff{e:yJsign} is replaced by
\[
  \frac{\partial y_1}{\partial J} \,>\, 0~, \quad
  \frac{\partial y_2}{\partial J} \,<\, 0~, \quad
  \frac{\partial y_1}{\partial J} + \frac{\partial y_2}{\partial J}
  \,=\, -\frac{\partial y_3}{\partial J} \,>\, 0~, \quad
  \hbox{if } J > 0~.
\]
In particular, if $J > 0$, we obtain
\[
  -\frac{\partial T}{\partial J} \,=\, (A_1-A_2)\frac{\partial y_1}
  {\partial J} + A_2 \Bigl(\frac{\partial y_1}{\partial J} +
  \frac{\partial y_2}{\partial J}\Bigr) \,>\, 0~,
\]
which proves ii).

The main difficulty is the proof of i), which requires a more
sophisticated argument. As $y_1 - y_3 < s(\phi) - y_3 < y_2 - y_3$
for all $\phi \in (0,\pi/2)$, it is clear that
\[
  \frac{c}{(y_2-y_3)^{3/2}} \,<\, A_1, A_2 \,<\,
  \frac{c}{(y_1-y_3)^{3/2}}~, \quad \hbox{where} \quad
  c \,=\, \frac{3\pi\sqrt{2}}4~.
\]
Using the upper bound for $A_1$ and the lower bound for
$A_2$ we obtain from \reff{e:Tder1}:
\begin{equation}\label{e:Tderr1}
  -\frac{\partial T}{\partial E} \,>\, \frac{c}{(y_1-y_3)^{3/2}}
  \,\frac{\partial y_1}{\partial E} + \frac{c}{(y_2-y_3)^{3/2}}
  \,\frac{\partial y_2}{\partial E}~.
\end{equation}
Replacing the derivatives $\partial y_1/\partial E$ and
$\partial y_2/\partial E$ with their expressions \reff{e:yiEJder},
we see that the right-hand side of \reff{e:Tderr1} is positive
if and only if
\[
  P_1 \,\eqdef\, y_1(y_2-y_3)^{3/2}(4E-3y_2^2-4y_2) +
  y_2(y_1-y_3)^{3/2}(4E-3y_1^2-4y_1) \,>\, 0~.
\]
Now, since $4E-3y_2^2-4y_2 < 0$ and $4E-3y_1^2-4y_1 > 0$, this
inequality is equivalent to
\[
  \frac{y_1|4E-3y_2^2-4y_2|}{y_2(4E-3y_1^2-4y_1)}
  \,<\, \Bigl(\frac{y_1-y_3}{y_2-y_3}\Bigr)^{3/2}~, \quad
  \hbox{where} \quad 0 < \frac{y_1-y_3}{y_2-y_3} < 1~.
\]
Clearly, a stronger inequality is obtained if we replace the
exponent $3/2$ by $2$ in the right-hand side. Thus it is sufficient
to show that $P_2 > 0$, where
\[
  P_2 \,=\, y_1(y_2-y_3)^2(4E-3y_2^2-4y_2) +
  y_2(y_1-y_3)^2(4E-3y_1^2-4y_1) \,>\, 0~.
\]
To do that, we recall that $y_1, y_2, y_3$ are the roots of the
cubic equation $y^3 +2y^2 - 4Ey + 2J^2 = 0$. In particular, we
have $y_1 + y_2 + y_3 = -2$ and $y_1y_2y_3 = -2J^2$. These
relations allow to eliminate the variables $E$ and $y_3$ from the
expression of $P_2$, which thus becomes a function of $y_1, y_2$
only. It is then convenient to set $y_1 = y-z$ and $y_2 =
y+z$, where $0 < z \le y$. After a straightforward algebra,
we obtain the final expression
\[
  P_2 \,=\, 8z^2 \Bigl((2+3y)^2 + z^2(3+4y)\Bigr)~,
\]
which shows that $P_2 > 0$. This concludes the proof. \QED

Using Proposition~\ref{th:monot1}, it is now easy to verify
that the Hamiltonian system associated to \reff{e:sgl1} is
non-degenerate, i.e. the determinant $\Delta(J,E)$ defined
in \reff{e:Deltadef} is always nonzero.

\begin{Proposition}\label{th:KAM1}
For all $(J,E) \in D$ we have $\Delta(J,E) > 0$.
\end{Proposition}

\noindent{\bf Proof.} Differentiating the renormalized phase
$\Psi$ with respect to $J$ and $E$, we obtain as in
\reff{e:Psider}
\begin{equation}\label{e:Psider1}
  \frac{\partial \Psi}{\partial E} \,=\, B_1 J\frac{\partial y_1}{\partial E}
  + B_2 J\frac{\partial y_2}{\partial E}~, \quad
  \frac{\partial \Psi}{\partial J} \,=\, B_1 J\frac{\partial y_1}{\partial J}
  + B_2 J\frac{\partial y_2}{\partial J} - B_3~, \quad
\end{equation}
where
\begin{eqnarray*}
  {B_1 \atop B_2} &=&
  \sqrt{2} \int_0^{\pi/2} \biggl\{\frac{1}{(-y_3)^{3/2}(s(\phi)-y_3)^{1/2}}
  \,\frac{2\sqrt{-y_3}+\sqrt{s(\phi)-y_3}}{(\sqrt{-y_3}
  +\sqrt{s(\phi)-y_3})^2}\,+ \\
  && \qquad \qquad \frac{1}{(-y_3)^{1/2}(s(\phi)-y_3)^{3/2}}
  \,\frac{\sqrt{-y_3}+2\sqrt{s(\phi)-y_3}}{(\sqrt{-y_3}+
  \sqrt{s(\phi)-y_3})^2}
  \,\Bigl(1+{\cos^2(\phi) \atop \sin^2(\phi)}\Bigr)\biggr\}\dd\phi~,
  \\[1mm]
  B_3 &=& 2\sqrt{2} \int_0^{\pi/2} \frac{1}{\sqrt{-y_3}
  \,\sqrt{s(\phi)-y_3}\,(\sqrt{-y_3} + \sqrt{s(\phi)-y_3})}\dd\phi
  \,=\, -\frac{\Psi}{J}~.
\end{eqnarray*}
Thus replacing \reff{e:Tder1}, \reff{e:Psider1} into \reff{e:Deltadef}
and proceeding as in the proof of Proposition~\ref{th:KAM}, we
obtain
\[
  \Delta(J,E) \,=\, (A_2 B_1 - A_1 B_2)(y_2-y_1)\frac{\partial y_1}
  {\partial J}\frac{\partial y_2}{\partial J} - B_3
  \frac{\partial T}{\partial E}~.
\]
As $\partial T/\partial E < 0$ by Proposition~\ref{th:monot1},
it is thus sufficient to verify that $\Delta_1 = A_2B_1 - A_1B_2 < 0$.
This inequality can be established using the same arguments as
in the the defocusing case. Indeed, if we define (for $i = 1,2$)
$\bA_i$, $\bB_i$ by \reff{e:bAiBidef} and $\tilde \bB_i$ by
\reff{e:tildeBidef}, we find $\Delta_1 = 2(-y_3)^{-4}\Delta_2$
where
\[
  \Delta_2  \,=\, \bA_2\bB_1 - \bA_1\bB_2 \,=\,
  \bA_2\tilde\bB_1 - \bA_1\tilde\bB_2~.
\]
Now, we observe that $\phi \mapsto \sigma(\phi) =
(1-y_3^{-1}s(\phi))^{1/2}$ is {\em increasing} over $[0,\pi/2]$,
because $y_3 < 0$. Using Lemma~\ref{th:fg}, we deduce that
$\bA_1 > \bA_2 > 0$ and $\tilde \bB_2 > \tilde \bB_1 > 0$, hence
$\Delta_2 < 0$. \QED

As in Section~\ref{s:ex}, Propositions~\ref{th:monot1},
\ref{th:KAM1} allow to determine the range of values of
the period $T$ and the renormalized phase $\Psi$. We find that
$(T,\Psi) : D \to \hat D$ is a smooth diffeomorphism, where
\[
  \hat D \,=\, \Big\{(T,\Psi) \in \R^2 \,\Big|\, 0 < T < \pi\,,~
  |\Psi| < \hat \Psi(T)\Big\}~, \quad \hbox{where} \quad
  \hat \Psi(T) \,=\, \pi - \Bigl(\frac{T^2+2\pi^2}{3}\Bigr)^{1/2}~.
\]
The domains $D$ and $\hat D$ are represented in Fig.~5.

Now, we fix $(J,E) \in D$ and we study the stability of the
periodic wave $U_{J,E}(x,t) = \rme^{-\rmi t} W_{J,E}(x)$, where
$W_{J,E}$ is a solution of \reff{e:sgl1} satisfying \reff{e:EJ1}.
As in Section~\ref{s:os}, we set $W_{J,E}(x)  = \rme^{\rmi px}
Q_{J,E}(2kx)$, where $k, \ell$ are defined in \reff{e:kldef} and
$p = k+\ell$. The discussion follows exactly the same lines as in
the defocusing case, so we shall just mention the main differences.
The function $Q(z,t)$ defined in \reff{e:Qgen} satisfies the
evolution equation
\begin{equation}\label{e:Qeq1}
  \rmi Q_t + 4\rmi pk Q_z + 4k^2 Q_{zz} + (1-p^2)Q + |Q|^2Q \,=\, 0~,
\end{equation}
and the corresponding energy functional reads
\begin{equation}\label{e:EEdef1}
  \EE(Q) \,=\, \int_0^{2\pi} \Bigl(2k^2 |Q_z(z)|^2
   - \frac14|Q(z)|^4\Bigr)\dd z~.
\end{equation}
In particular, if we define $\EE_{J,E}$ by \reff{e:Eabdef},
the second variation becomes
\[
  H_{J,E} \,=\, \EE_{J,E}''(Q_{J,E}) \,=\, -4k^2 \partial_{zz}
  -4\rmi pk\partial_z - (1{-}p^2) - |Q_{J,E}|^2 - 2 Q_{J,E}
  \otimes Q_{J,E}~.
\]

As in the defocusing case, we rely on the result found for small
waves in \cite[Remark A.2]{GH1} and conclude that, when $J = 0$
and $E > 0$ is sufficiently small, the operator $H_{J,E}$ acting
on $X_0$ has exactly three eigenvalues $\{\lambda_1,\lambda_2,0\}$
in a neighborhood of the origin, where
\[
  \lambda_1(E) \,=\, -3E + \rmO(E^2)~, \quad
  \lambda_2(E) \,=\, E + \rmO(E^2)~, \quad
  \hbox{as } E \to 0~.
\]
The other eigenvalues of $H_{J,E}$ are positive and bounded away
from zero. On the other hand, we know that zero is an eigenvalue
of $H_{J,E}$ of multiplicity exactly two for all $(J,E) \in D$.
Indeed, this property was established in Section~\ref{s:os} by a
general argument which uses only the symmetries of \reff{e:sgl}
and not the particular form of the nonlinearity. By continuity, it
follows that $H_{J,E}$ has exactly one negative eigenvalue for all
$(J,E) \in D$, so that Proposition~\ref{th:Hneg} remains valid in
the present case.

Finally, we recall that the equilibrium $Q_{J,E}$ of \reff{e:Qeq1}
is a member of a two-parameter family of travelling and rotating
waves of the form \reff{e:travrot}. For sufficiently small
$(\omega,c)$, the wave profile is given by $Q_{J,E}^{\omega,c}(z) =
\lambda Q_{J'\!,E'}(z)$, where $\lambda$ is defined by
\reff{e:lamvdef} and $(J',E') \in D$ is the only point in
a neighborhood of $(J,E)$ such that \reff{e:omcdef} holds.
Let $\HH_{J,E}$ be the Hessian matrix of the function $d_{J,E}$
defined by \reff{e:ddef}, \reff{e:EJEomc}. In view of
Proposition~\ref{th:Hneg}, the general results of \cite{GSS2}
imply:

\begin{Proposition}\label{th:orbit1}
For all $(J,E) \in D$ such that $\det(\HH_{J,E}) < 0$, the
periodic wave $Q_{J,E}$ is a stable equilibrium of \reff{e:Qeq1}
in the sense of Proposition~\ref{th:Qorbit}.
\end{Proposition}

\figurewithtex 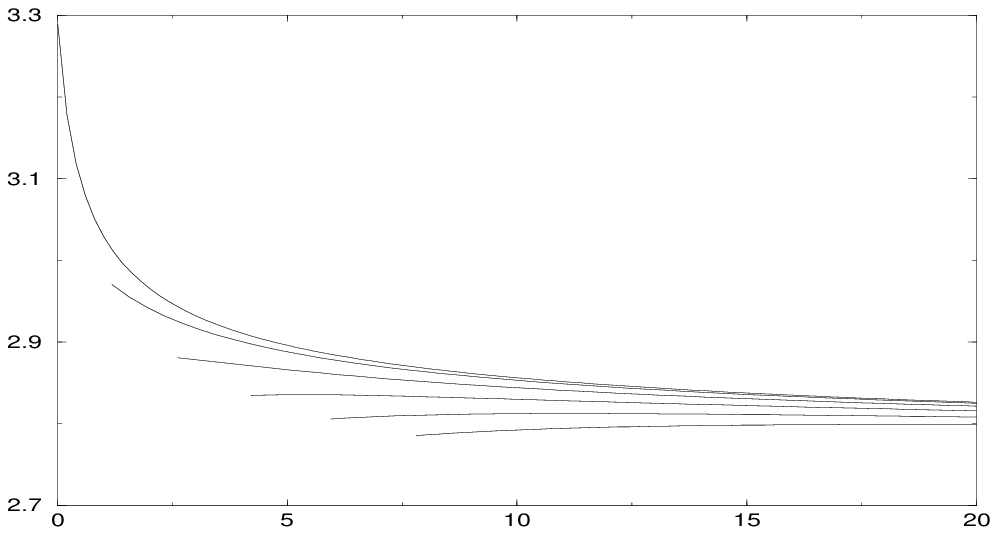 Fig6.tex 5.000 10.000
{\bf Fig.~6:} The quantity $-\det(\HH_{J,E})$ is represented
as a function of $E \in [E_-(J),20]$ for $J = 0,1,\dots,5$.
Similar curves are obtained for higher values of $J$, thus
indicating that $\det(\HH_{J,E})$ is always negative.\cr

Proposition~\ref{th:orbit1} is a conditional stability result,
since it applies under the assumption that $\det(\HH_{J,E}) < 0$.
This condition is satisfied at least for $(J,E)$ sufficiently
close to zero, because we know from \cite{GH1} that
\[
  \HH_{J,E} \,=\, \frac{\pi}3 \pmatrix{2 & 1 \cr 1 & -1}
  (\1 + \rmO(E))~, \quad \hbox{as} \quad (J,E) \to (0,0)
  \hbox{ in }D~.
\]
Moreover, the Hessian matrix $\HH_{J,E}$ is easy to evaluate
numerically for any $(J,E) \in D$, and its determinant appears to
be always negative (see Fig.~6). Thus we conjecture that
$\det(\HH_{J,E}) < 0$ for all $(J,E) \in D$. This property can
probably be established rigorously using similar arguments as in
the proof of Proposition~\ref{th:Hessian}, but the modifications
are not straightforward.

\subsection{Corotating waves ($\omega = -1$)}
\label{s:om2}

Finally we study the corotating waves of the focusing nonlinear
Schr\"odinger equation. Our starting point is the stationary equation
\begin{equation}\label{e:sgl2}
  W_{xx}(x) - W(x) + |W(x)|^2 W(x) \,=\, 0~, \quad x \in \R~.
\end{equation}
The invariants of this Hamiltonian system have the following
expressions:
\begin{equation}\label{e:EJ2}
  J \,=\, \Im(\overline{W}W_x)~, \qquad
  E \,=\, \frac12|W_x|^2 - \frac12|W|^2 + \frac14|W|^4~.
\end{equation}
It is convenient to use the parametrization $J = q(1+q^2)$,
where $q \in \R$. If $J \neq 0$, the effective potential
$V_J(r) = J^2/(2r^2) - r^2/2 + r^4/4$ has a unique critical
point at $r = r_q = \sqrt{1+q^2}$, where $V_J$ attains its
global minimum:
\[
  E_-(J) \,=\, V_J(\sqrt{1+q^2}) \,=\, \frac14 (q^2+1)(3q^2-1)~.
\]
In that case, Eq.\reff{e:sgl2} has quasi-periodic solutions
for all $E > E_-(J)$. If $J = 0$, the double-well potential
$V_0(r)$ has two minima at $r = \pm 1$ and a local maximum at
$r = 0$. It follows that \reff{e:sgl2} has (real) periodic
solutions if $-1/4 < E < 0$ (dnoidal waves) and if $E > 0$ 
(cnoidal waves). Summarizing, the parameter domain where 
quasi-periodic solutions of \reff{e:sgl2} exist is
\[
  D \,=\, \{ (J,E) \in \R^2 \,|\, E > E_-(J)\} \setminus
  \{(0,0)\}~.
\]
In the exceptional case $(J,E) = (0,0)$, Eq.~\reff{e:sgl2}
has the pulse-like solution $W(x) = \sqrt{2}/\cosh(x)$ which
corresponds to the solitary wave of the focusing NLS equation.

\figurewithtex 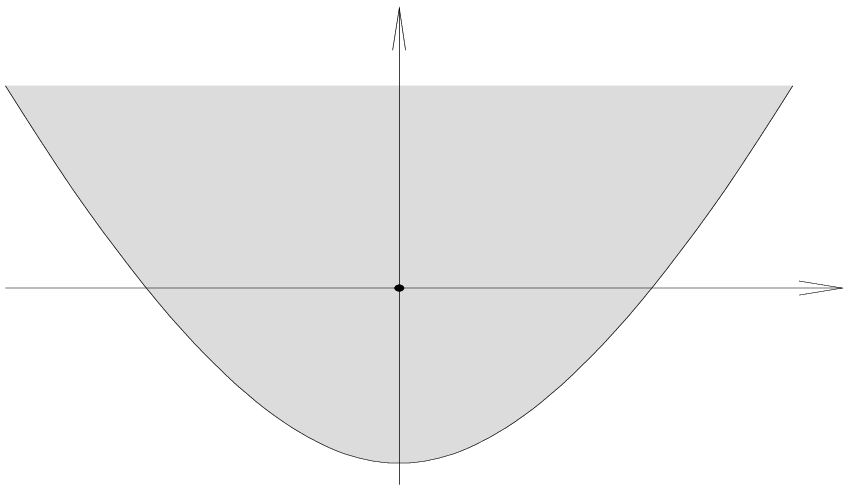 Fig7.tex 5.00 9.00
{\bf Fig.~7:} Existence domain for the corotating waves of
the focusing NLS equation. The origin $(J,E)=(0,0)$ corresponds
to the solitary wave. The half-line $\Gamma_+$ (cnoidal waves)
and the segment $\Gamma_-$ (dnoidal waves) are the discontinuity 
lines of the phases $\Phi$ and $\Psi$, respectively.\cr

If $(J,E) \in D$, the period $T$, the phase increment $\Phi$,
and the renormalized phase $\Psi$ of the quasi-periodic solutions
are given by the formulas \reff{e:TPhiexp1}, \reff{e:Psiexp1},
where $s(\phi) = y_1 \cos^2(\phi) + y_2 \sin^2(\phi)$ and
$y_3 \le 0 \le y_1 < y_2$ are the roots of the cubic polynomial
$P(y) = -y^3 +2y^2 + 4Ey - 2J^2$. It is important to realize
that the phase $\Phi$ is continuous on $D$, except on the
half-line
\[
  \Gamma_+ \,=\, \{(0,E) \in \R^2 \,|\, E > 0\} \,=\,
  \{(J,E) \in D \,|\, y_1 = 0\}~,
\]
(see Fig.~7). Similarly, the renormalized phase $\Psi$ is
continuous on $D$ except for a $2\pi$-jump on the line segment
\[
  \Gamma_- \,=\, \{(0,E) \in \R^2 \,|\, -1/4 < E < 0\} \,=\,
  \{(J,E) \in D \,|\, y_3 = 0\}~.
\]
We shall thus use either $\Phi$ or $\Psi$ depending on the
parameter region under consideration.

Unlike in the previous cases, the period $T$ is no longer a
monotone function of the energy $E$. This is intuitively
clear, as we expect that $T \to +\infty$ as $(J,E) \to (0,0)$.
In fact, one can prove:

\begin{Lemma}\label{th:notmonot}
The period $T(J,E)$ satisfies
\[
  \frac{\partial T}{\partial E}(0,E) > 0 \hbox{ for }
  -1/4 < E < 0~, \quad \hbox{and} \quad
  \frac{\partial T}{\partial E}(0,E) < 0 \hbox{ for }
  E > 0~.
\]
\end{Lemma}

\noindent{\bf Proof.} As in \reff{e:Tder1}, we have
\[
  \frac{\partial T}{\partial E} \,=\, -A_1 \frac{\partial y_1}
  {\partial E} -A_2 \frac{\partial y_2}{\partial E}~,
\]
where the coefficients $A_1$, $A_2$ are defined in \reff{e:Aiexp1}.
The only difference with the previous case is that the quantities
$y_1,y_2,y_3$ which appear in \reff{e:Aiexp1} are the roots of
a different polynomial. If $J = 0$ and $E > 0$, then $y_1 = 0$,
$y_2 = 1+\sqrt{1+4E}$, $y_3 = 1-\sqrt{1+4E}$, hence
\[
  \frac{\partial T}{\partial E}(0,E) \,=\,-A_2 \,\frac{2}{\sqrt{1+4E}}
  \,<\, 0~.
\]
If $J = 0$ and $-1/4 < E < 0$, then $y_1 = 1-\sqrt{1+4E}$,
$y_2 = 1+\sqrt{1+4E}$, $y_3 = 0$, so that
\[
  \frac{\partial T}{\partial E}(0,E) \,=\,\frac{2}{\sqrt{1+4E}}\,
  (A_1 - A_2) \,>\, 0~,
\]
because $A_1 > A_2$ by \reff{e:A1A21}. \QED

On the other hand, the period $T$ is still a monotone function
of $J$ if $J > 0$. As in Proposition~\ref{th:monot1}, one can
prove that
\[
   \frac{\partial \Psi}{\partial E}(J,E) \,=\,
   \frac{\partial \Phi}{\partial E}(J,E) \,=\,
   -\frac{\partial T}{\partial J}(J,E) \,>\, 0~, \quad
   \hbox{if } J > 0~.
\]
The important quantity is again
\begin{equation}\label{e:Deltadef2}
  \Delta(J,E) \,=\, \Det\left(\begin{array}{cc}
  \frac{\partial T}{\partial E} &
  \frac{\partial \Phi}{\partial E} \\[1mm]
  \frac{\partial T}{\partial J} &
  \frac{\partial \Phi}{\partial J}
  \end{array}\right)(J,E) \,=\,
  \Det\left(\begin{array}{cc}
  \frac{\partial T}{\partial E} &
  \frac{\partial \Psi}{\partial E} \\[1mm]
  \frac{\partial T}{\partial J} &
  \frac{\partial \Psi}{\partial J}
  \end{array}\right)(J,E)~.
\end{equation}
Here the first determinant in the right-hand side is meaningful
if $(J,E) \notin \Gamma_+$, and the second one if $(J,E) \notin
\Gamma_-$.

\begin{Proposition}\label{th:KAM2}
For all $(J,E) \in D$ we have $\Delta(J,E) > 0$.
\end{Proposition}

\noindent{\bf Proof.}
Fix $(J,E) \in D$. If $(\partial T/\partial E)(J,E) < 0$, we
know from Lemma~\ref{th:notmonot} that $(J,E) \notin \Gamma_-$.
Thus, using the second determinant in \reff{e:Deltadef2} and
proceeding exactly as in Proposition~\ref{th:KAM1}, we obtain
$\Delta(J,E) > 0$.

We now assume that $(\partial T/\partial E)(J,E) \ge 0$,
so that $(J,E) \notin \Gamma_+$ by Lemma~\ref{th:notmonot}.
Differentiating the expression of $\Phi$ in \reff{e:TPhiexp1}
with respect to $E$ and $J$, we obtain
\begin{equation}\label{e:Phider2}
  \frac{\partial \Phi}{\partial E} \,=\,
  -\BB_1 J \frac{\partial y_1}{\partial E}
  -\BB_2 J\frac{\partial y_2}{\partial E}~, \quad
  \frac{\partial \Phi}{\partial J} \,=\,
  -\BB_1 J\frac{\partial y_1}{\partial J}
  -\BB_2 J\frac{\partial y_2}{\partial J} + \BB_3~,
\end{equation}
where
\begin{eqnarray*}
  \BB_1 &=& \sqrt{2} \int_0^{\pi/2}\Bigl(
  \frac{2\cos^2(\phi)}{s(\phi)^2\sqrt{s(\phi)-y_3}} +
  \frac{1+\cos^2(\phi)}{s(\phi)(s(\phi)-y_3)^{3/2}}
  \Bigr)\dd \phi ~=~ \BB_{11} + \BB_{12}~, \\
  \BB_2 &=& \sqrt{2} \int_0^{\pi/2}\Bigl(
  \frac{2\sin^2(\phi)}{s(\phi)^2\sqrt{s(\phi)-y_3}} +
  \frac{1+\sin^2(\phi)}{s(\phi)(s(\phi)-y_3)^{3/2}}
  \Bigr)\dd \phi ~=~ \BB_{21} + \BB_{22}~, \\
  \BB_3 &=& 2\sqrt{2} \int_0^{\pi/2} \frac{1}{s(\phi)
  \sqrt{s(\phi)-y_3}}\dd\phi ~=~ \frac{\Phi}{J}~.
\end{eqnarray*}
Then, using \reff{e:Tder1}, \reff{e:Phider2}, and the first determinant
in \reff{e:Deltadef2}, we find as in Proposition~\ref{th:KAM}:
\[
  \Delta(J,E) \,=\, (A_1 \BB_2 - A_2 \BB_1)(y_2-y_1)
  \frac{\partial y_1}{\partial J}
  \frac{\partial y_2}{\partial J} + \BB_3
  \frac{\partial T}{\partial E}~.
\]
As $(\partial T/\partial E) \ge 0$ by assumption, it is sufficient
to verify that $\Delta_1 \eqdef A_2 \BB_1 - A_1 \BB_2 > 0$.

To do that, we observe that $\Delta_1 = \Delta_{11} + \Delta_{12}$
where
\[
  \Delta_{11} \,=\, A_2 \BB_{11} - A_1 \BB_{21}~, \quad
  \Delta_{12} \,=\, A_2 \BB_{12} - A_1 \BB_{22}~,
\]
and we prove separately that $\Delta_{11} > 0$ and $\Delta_{12} > 0$.
Both inequalities are easy consequences of Lemma~\ref{th:fg}.
Indeed, $\Delta_{11} > 0$ is equivalent to $(A_1{-}A_2)/(A_1{+}A_2)
< (\BB_{11}{-}\BB_{21})/(\BB_{11}{+}\BB_{21})$, or explicitly
\[
  \frac{\displaystyle{\int_0^{\pi/2}\frac{\cos^2(\phi)-\sin^2(\phi)}
  {(s(\phi)-y_3)^{3/2}} \dd \phi}}{\displaystyle{\int_0^{\pi/2}\frac{3}
  {(s(\phi)-y_3)^{3/2}} \dd \phi}} \quad < \quad
  \frac{\displaystyle{\int_0^{\pi/2}\frac{\cos^2(\phi)-\sin^2(\phi)}
  {s(\phi)^2(s(\phi)-y_3)^{1/2}} \dd \phi}}{\displaystyle{\int_0^{\pi/2}
  \frac{1}{s(\phi)^2(s(\phi)-y_3)^{1/2}}\dd \phi}}~.
\]
A stronger inequality (without the factor $3$ in the denominator
of the left-hand side) is obtained from Lemma~\ref{th:fg} by
choosing $I = [0,\pi/2]$, $f(\phi) = \cos^2(\phi) - \sin^2(\phi)$,
$g(\phi) = (s(\phi){-}y_3)/(s(\phi))^2$ and
\[
  \dd \mu \,=\, \frac{1}{\NN}\,\frac{\dd \phi}{(s(\phi)-y_3)^{3/2}}~,
  \quad \hbox{where}\quad \NN \,=\, \int_0^{\pi/2}
  \frac{\dd \phi}{(s(\phi)-y_3)^{3/2}}~.
\]
Indeed $\int_I fg\dd\mu > (\int_I f\dd\mu)(\int_I g\dd\mu)$
because $f,g$ are strictly decreasing over $I$. Thus $\Delta_{11}>0$,
and the same argument with $g(\phi) = 1/s(\phi)$ shows that
$\Delta_{12} > 0$. \QED

To conclude this section, we fix $(J,E) \in D$ and we study the
stability of the periodic wave $U_{J,E}(x,t) = \rme^{\rmi t}
W_{J,E}(x)$, where $W_{J,E} = \rme^{\rmi px} Q_{J,E}(2kx)$ is a
solution of \reff{e:sgl2} satisfying \reff{e:EJ2}. Setting
$U(x,t) = \rme^{\rmi (px+t)}Q(2kx,t)$, we obtain from \reff{e:nls}
with $\gamma = 1$ the evolution equation
\begin{equation}\label{e:Qeq2}
  \rmi Q_t + 4\rmi pk Q_z + 4k^2 Q_{zz} - (1+p^2)Q + |Q|^2Q \,=\, 0~.
\end{equation}
The wave profile $Q_{J,E}$ is a critical point of the modified
energy
\[
  \EE_{J,E}(Q) \,=\, \EE(Q) + (1+p^2)N(Q) - 4p k M(Q)~,
\]
where $\EE$ is defined in \reff{e:EEdef1} and $N,M$ in
\reff{e:NMEdef}. The second variation of $\EE_{J,E}$ at $Q_{J,E}$
is
\[
  H_{J,E} \,=\, \EE_{J,E}''(Q_{J,E}) \,=\, -4k^2 \partial_{zz}
  -4\rmi pk\partial_z + (1{+}p^2) - |Q_{J,E}|^2 - 2 Q_{J,E}
  \otimes Q_{J,E}~.
\]
Again, one can prove that Proposition~\ref{th:Hneg} still holds in the
present case. The fact that zero is always a double eigenvalue of
$H_{J,E}$ is established as in Section~\ref{s:os}, and a direct
calculation for small amplitude periodic waves (in a neighborhood of a
plane wave) shows that $H_{J,E}$ has exactly one negative
eigenvalue.

As in the previous cases, the equilibrium $Q_{J,E}$ of \reff{e:Qeq1}
is a member of a two-parameter family of travelling and rotating
waves of the form \reff{e:travrot}. For sufficiently small
$(\omega,c)$, the wave profile is given by $Q_{J,E}^{\omega,c}(z) =
\lambda Q_{J'\!,E'}(z)$, where $\lambda$ is defined by
\reff{e:lamvdef} and $(J',E') \in D$ is the only point in
a neighborhood of $(J,E)$ such that
\[
  \omega \,=\, (1+p^2) - \lambda^2(1+p'^2)~, \quad
  c \,=\, 4 \lambda^2 k'p' - 4 kp~.
\]
If $\HH_{J,E}$ denotes the Hessian matrix of the function $d_{J,E}$
defined by \reff{e:ddef}, \reff{e:EJEomc}, the results of \cite{GSS2}
imply:

\begin{Proposition}\label{th:orbit2}
For all $(J,E) \in D$ such that $\det(\HH_{J,E}) < 0$, the
periodic wave $Q_{J,E}$ is a stable equilibrium of \reff{e:Qeq2}
in the sense of Proposition~\ref{th:Qorbit}.
\end{Proposition}

As in Section~\ref{s:om1}, we conjecture that $\det(\HH_{J,E}) < 0$
for all $(J,E) \in D$. This inequality is true at least for small
amplitude periodic waves (in a neighborhood of a plane wave), and
numerical calculations indicate that it remains valid over the whole
parameter domain $D$. In the particular case where $(J,E) \in
\Gamma_-$ (dnoidal waves), the orbital stability with respect to
periodic perturbations has been established in \cite{Ang}. Remark
that Proposition~\ref{th:orbit2} does apply in the case where $(J,E)
\in \Gamma_+$ (cnoidal waves), which is not covered by the results of
\cite{Ang}, but as is explained in the introduction this is because we
use in fact a more restricted class of perturbations.

\end{document}